\documentclass{article}

\usepackage{vmargin,epsfig}
\usepackage[english,francais]{babel}
\usepackage[latin1]{inputenc}
\usepackage[tbtags]{amsmath}
\usepackage{amsthm,amssymb}
\usepackage{pifont}
\usepackage{bbm}
\usepackage[all]{xy}

\newtheorem{theo}{Theorème}[section]
\newtheorem{lemme}[theo]{Lemme}
\newtheorem{prop}[theo]{Proposition}

\theoremstyle{definition}
\newtheorem{deftn}[theo]{Définition}
\newtheorem{conj}[theo]{Conjecture}

\def\leq{\leqslant}
\def\geq{\geqslant}

\def\Z{\mathbb{Z}}
\def\Q{\mathbb{Q}}

\def\F{\mathbb{F}}
\def\O{\mathcal{O}}

\def\G{\mathcal{G}}

\def\pa#1{\left(#1\right)}
\def\cro#1{\left[ #1 \right]}

\def\brac#1{\left< #1 \right>}

\def\epsilon{\varepsilon}

\def\calA{\mathcal{A}}
\def\calB{\mathcal{B}}
\def\calC{\mathcal{C}}
\def\calE{\mathcal{E}}
\def\calF{\mathcal{F}}
\def\calM{\mathcal{M}}
\def\calN{\mathcal{N}}

\def\id{\text{\rm id}}

\def\hom{\text{\rm Hom}}
\def\Fil{\text{\rm Fil}}
\def\Gal{\text{\rm Gal}}

\def\Rep{\text{\rm Rep}}

\def\GL{\text{\rm GL}}
\def\Frob{\text{\rm Frob}}
\def\det{\text{\rm det}}
\def\Sym{\text{\rm Sym}}

\def\Gr{\text{\rm Gr}}
\def\EGrdd{E\text{\rm -Gr}_{\text{\rm dd}}}
\def\EBrModdd{\text{\rm E-BrMod}_{\text{\rm dd}}}

\def\Vcris{V_{\text{\rm cris}}^\star}
\def\Tcris{T_{\text{\rm cris}}^\star}
\def\Tst{T_{\text{\rm st}}^\star}
\def\f{\text{\rm fil}}

\def\det{\text{\rm det}}
\def\Spec{\text{\rm Spec}}

\def\cris{\text{\rm cris}}

\def\1{\mathbbm{1}}

\title{Schémas en groupes et poids de Diamond-Serre}
\author{Xavier Caruso}
\date{Avril 2007}

\begin{document}

\maketitle

\tableofcontents

\bigskip

\noindent \hrulefill

\bigskip

Cette note fait suite à un exposé que j'ai donné pour la deuxième 
partie du groupe de travail sur \emph{les généralisations de la 
conjecture de modularité de Serre} organisé par Christophe Breuil, Guy 
Henniart, Ariane Mézard et Rachel Ollivier au printemps 2007. Le thème
de l'exposé est la proposition 3.3.1 de \cite{gee} : il s'agit d'abord
de présenter le matériel de théorie de Hodge $p$-adique nécessaire pour
comprendre l'énoncé de cette proposition, puis d'en donner une 
démonstration.

Nous avons jugé utile d'écrire ce texte car la rédaction de \cite{gee}
est souvent lacunaire ou imprécise et, surtout, contient une erreur :
telle qu'énoncée dans \cite{gee}, la proposition 3.3.1 est fausse, comme
l'auteur de cette note s'en est rendu compte lors de la préparation de
l'exposé. Il est en fait assez délicat de comprendre ce qu'il faut
modifier pour rendre l'énoncé juste, et c'est Gee lui-même qui a proposé
la bonne correction qui consiste à changer la définition de classe $J$
(définition \ref{def:classeJ}) en introduisant un décalage d'indice.
Cela fait, il n'est pas non plus immédiat de saisir comment les
arguments de \cite{gee} s'adaptent à ce nouvel énoncé. Cette partie du
travail a été accomplie par l'auteur de cette note.

Nous présentons dans ce papier la version modifiée de la proposition 3.3.1 
accompagnée de sa démonstration (correcte, nous l'espérons). Notons pour 
finir que Gee affirme que c'est bien cette version modifiée qu'il utilise
par la suite dans son article ; ainsi, la discussion précédente ne remet
\emph{a priori} pas en cause les résultats de \cite{gee}. Nous espérons que 
ce texte pourra être utile à tous ceux qui veulent apprendre la théorie.

\section{Rappel du contexte}

La proposition 3.3.1 de \cite{gee} est un résultat de théorie de Hodge
$p$-adique pure, et peut tout à fait être énoncée (et prouvée) sans 
aucune référence à la conjecture de modularité pour les corps totalement 
réels énoncée par Buzzard, Diamond et Jarvis dans \cite{bdj}. 
Toutefois, en procédant ainsi, il deviendrait difficile de comprendre
l'intérêt de cet énoncé, et c'est pourquoi nous préférons consacrer ce
premier chapitre à un bref rappel du contexte\footnote{Qui, lors du 
groupe de travail, a fait l'objet d'un exposé complet d'Ariane Mézard.},
ce qui permettra par là-même de fixer les notations\footnote{Certaines
notations de \cite{gee} diffèrent de celles de \cite{bdj}. Nous 
employons naturellement dans cette note celles de \cite{gee}.}.

\subsection{La conjecture de Buzzard, Diamond et Jarvis}

Fixons $p$ un nombre premier et $F$ un corps de nombres totalement réel 
dans lequel $p$ est non ramifié. Notons $\O_F$ l'anneau des entiers de
$F$ et $G_F$ son groupe de Galois absolu. Fixons également $\bar 
\F_p$ une clôture algébrique du corps fini à $p$ éléments $\F_p$.
Soit $\rho : G_F \to \GL_2(\bar \F_p)$ une représentation \emph{que l'on 
supposera toujours continue, irréductible et totalement impaire.} La 
généralisation naturelle de la conjecture de modularité de Serre est la 
suivante :

\begin{conj}
\label{conj:1}
Avec les notations précédentes, $\rho$ est modulaire, c'est-à-dire 
qu'elle est isomorphe à la réduction modulo $p$ d'une représentation 
associée à une forme modulaire de Hilbert.
\end{conj}

Le travail de Buzzard, Diamond et Jarvis a consisté à produire une 
version raffinée de la conjecture précédente, dans laquelle on prédit en 
outre le poids de la forme modulaire. Pour énoncer cette forme plus 
précise, on commence par définir la notion de poids utilisée dans ce 
contexte :

\begin{deftn}
Un \emph{poids de Diamond-Serre} est une (classe d'isomorphisme de) 
représentation(s) irréductible(s) de $\GL_2(\O_F/p)$ à coefficients dans 
$\bar \F_p$.
\end{deftn}

\noindent
Buzzard, Diamond et Jarvis définissent\footnote{On donnera dans la 
suite de l'exposé quelques morceaux de ces définitions. Pour une 
définition complète, on renvoie le lecteur aux sections 2 et 3 de 
\cite{bdj}.} d'une part la notion d'être modulaire de poids $\sigma$ 
(pour une représentation vérifiant les mêmes hypothèses que $\rho$) et 
d'autre part un certain ensemble de poids (de Diamond-Serre) dépendant 
de $\rho$ et noté $W(\rho)$. Ils émettent alors la conjecture :

\begin{conj}
\label{conj:2}
La représentation $\rho$ est modulaire de poids $\sigma$ si, et 
seulement si $\sigma \in W(\rho)$.
\end{conj}

\noindent
Ils démontrent de plus que la modularité de poids $\sigma$ (pour un
$\sigma$) implique la modularité (tout court) et que l'ensemble
$W(\rho)$ n'est jamais vide. Ainsi la conjecture \ref{conj:2}
implique-t-elle la conjecture \ref{conj:1}. Le résultat principal de
\cite{gee} en est une réciproque partielle : sous certaines hypothèses
supplémentaires sur $\rho$, la conjecture \ref{conj:1} implique
l'équivalence de la conjecture \ref{conj:2} pour certains poids
$\sigma$, dits réguliers (voir définition \ref{def:regulier}).

\subsection{Énoncé (vague) de la proposition 3.3.1}

Nous commençons par donner quelques indications (celles qui seront
nécessaires pour l'exposé) sur la définition de $W(\rho)$. Par le lemme
chinois, le groupe $\GL_2(\O_F/p)$ se décompose comme le produit des
$\GL_2(k_v)$ où $v$ parcourt l'ensemble des idéaux premiers au-dessus de
$p$ et où $k_v$ désigne le corps résiduel associé. Tout poids $\sigma$
se décompose comme un produit tensoriel de représentations irréductibles
de $\GL_2(k_v)$. La définition de $W(\rho)$ est locale en ce sens que
l'on définit d'abord des ensembles $W_v(\rho)$ de représentations
irréductibles de $\GL_2(k_v)$ et que l'on convient ensuite que $\sigma
\in W(\rho)$ si, et seulement si ses facteurs locaux appartiennent tous
aux ensembles $W_v(\rho)$ respectifs. De plus, $W_v (\rho)$ ne dépend
que de la restriction de $\rho$ au groupe de décomposition en $v$.

\paragraph{Notations}

\emph{À partir de maintenant, et jusqu'à la fin de cette note}, nous
fixons une place $v$ au-dessus de $p$. Ainsi nous pourrons nous
affranchir des indices $v$. Par exemple, nous notons plus simplement $k$
à la place de $k_v$ : c'est une extension finie de $\F_p$ dont nous
désignons par la lettre $r$ le degré. Introduisons encore quelques
notations. Soit $K_0$ le complété de $F$ en la place $v$, c'est une
extension finie non ramifiée de $\Q_p$ de corps résiduel $k$. Soit $K =
K_0(\pi)$ où $\pi$ est une racine $(p^r-1)$-ième de $(-p)$, c'est une
extension galoisienne\footnote{En effet, $k \simeq \F_{p^r}$ contient
les racines $(p^r-1)$-ième de l'unité et donc $K_0$ aussi.} totalement
ramifiée de $K_0$. On note $G$ le groupe de Galois de cette extension,
et $e$ son degré, \emph{i.e.} $e = p^r -1$. Soit $\O_{K_0}$ (resp
$\O_K$) l'anneau des entiers de $K_0$ (resp. de $K$). Soient $G_{K_0}$
(resp. $G_K$) le groupe de Galois absolu de $K_0$ (resp. de $K$) et
$I_{K_0}$ (resp. $I_K$) son sous-groupe d'inertie. On remarque que
$\rho$ induit par restriction une représentation de $G_{K_0}$ (qui peut
être, elle, non irréductible). Dans la suite, on notera simplement
$\rho$ pour $\rho_{|G_{K_0}}$ et on ne considèrera plus la
représentation globale.

Définissons le caractère galoisien\footnote{Il n'est pas paru clair à 
l'auteur de cette note, si dans \cite{gee}, le caractère $\omega$ et
les $\omega_i$ qui en découlent sont ceux définis ici ou leurs 
inverses : Gee semble dire qu'il s'agit des inverses, mais certaines
utilisations qu'il en fait laisse croire le contraire. Nous choisissons
cette définition ici, car elle est plus adaptée pour notre propos, et
en plus en accord avec l'exposé précédent et les notations de 
\cite{bdj}.} $\omega : G_{K_0} \to k^\star$, $g \mapsto 
\frac{g(\pi)}{\pi}$ où le calcul est effectué dans $K$ puis 
réduit modulo $\pi$. Nous considérerons souvent $\omega$ comme un 
caractère de $I_{K_0}$ (par restriction) ou de $G$ (par passage au 
quotient). Soit $S = \Z/r\Z$. Fixons un plongement $\tau_0 : k \to \bar 
\F_p$ et pour tout $i \in S$, notons $\tau_i = \tau_0 \circ \Frob^{-i}$
où $\Frob : x \mapsto x^p$ est le Frobenius arithmétique.
Pour finir, posons $\omega_i = \tau_i \circ \omega$ pour tout $i \in 
S$. On a les relations $\tau_{i+1}^p = \tau_i$ et $\omega_{i+1}^p = 
\omega_i$.

\paragraph{Définition de l'ensemble $W_v(\rho)$}

On rappelle que les représentations irréductibles de $\GL_2(k)$
s'écrivent comme suit :
$$\bigotimes_{i \in S} \det^{a_i} \: \Sym^{b_i-1} k^2 \otimes_{\tau_i} 
\bar \F_p$$
où les $a_i$ et les $b_i$ sont des entiers avec $1 \leq b_i \leq p$. On 
peut également choisir les $a_i$ dans l'intervalle $[0, p-1]$ de telle 
sorte que tous ne soient pas égaux à $p-1$. Avec cette condition 
supplémentaire, les représentations écrites précédemment sont deux à 
deux non isomorphes.

\begin{deftn}
\label{def:regulier}
Un poids de Diamond-Serre est dit \emph{régulier en $v$} si, sur la 
description précédente, on a $2 \leq b_i \leq p-2$ pour tout $i \in S$.

Il est dit \emph{régulier} s'il est régulier en toutes les places $v$.
\end{deftn}

On suppose à partir de maintenant de $\rho$ est réductible de la forme
$\rho \simeq \pa{\begin{array}{cc} \psi' & \star \\ 0 & 
\psi'' \end{array} }$ pour des caractères $\psi'$ et $\psi''$ de 
$G_{K_0}$.

\begin{deftn}
\label{def:wrho}
L'ensemble $W(\rho)$ est constitué des représentations irréductibles 
$\sigma$ de $\GL_2(k)$ (dont on note $a_i$ et $b_i$ les entiers 
associés), telles qu'il existe $J$ un sous-ensemble de $S$ pour lequel 
les deux conditions suivantes sont satisfaites :

\medskip

\begin{itemize}
\item[\it Condition {\bf A}] :
$\displaystyle \psi'_{|I_{K_0}} \simeq \prod_{i \in S} \omega_i^{a_i} 
\prod_{i \in J} \omega_i^{b_i}$ et $\displaystyle \psi''_{|I_{K_0}} 
\simeq \prod_{i \in S} \omega_i^{a_i} \prod_{i \not\in J} 
\omega_i^{b_i}$
\item[\it Condition {\bf B}] :
\og $\rho$ admet un relevé cristallin d'un certain type. \fg
\end{itemize}
\end{deftn}

\medskip

\noindent
Pour des raisons de commodité qui apparaitront plus clairement vers la
fin de cette note, on introduit $\1_J$ et $\bar \1_J$ les fonctions 
indicatrices respectives des ensembles $J$ et $S \backslash J$. Ainsi, 
par exemple, la condition {\bf A} se réécrit sous la forme suivante :
$$ \psi'_{|I_{K_0}} \simeq \prod_{i \in S} \omega_i^{a_i + b_i \1_j(i)}
\quad \text{et} \quad \psi''_{|I_{K_0}} \simeq \prod_{i \in S} 
\omega_i^{a_i + b_i \bar \1_j(i)}.$$

\paragraph{But de l'exposé}

L'objet de la proposition 3.3.1 de \cite{gee} est de donner, dans le cas 
des poids réguliers, une expression de la condition {\bf B} en termes 
de schémas en groupes. Schématiquement, cette proposition prend la forme 
qui suit :

\begin{prop}
\label{prop:vague}
On considère un poids $\sigma$ pour lequel $a_i = 0$ et $2 \leq b_i \leq 
p-2$ pour tout $i$. On fixe un sous-ensemble $J \subset S$ et on suppose 
que la condition {\bf A} est satisfaite. Alors la condition {\bf B} est
équivalente à la :

\medskip

\begin{itemize}
\item[\it Condition {\bf B'}] :
\og $\rho_{|G_K}$ est la fibre générique d'un schéma en groupes sur 
$\O_K$ muni d'une donnée de descente de $K$ à $K_0$, d'un certain type. 
\fg
\end{itemize}
\end{prop}

\noindent
Évidemment, il reste à préciser le sens des deux occurrences de 
l'expression \og un certain type \fg. Pour cela, on a besoin 
d'introduire certaines notions de théorie de Hodge $p$-adique, et
notamment la théorie de Breuil. C'est l'objet de la section suivante.
Les sections \ref{sec:Bp} et \ref{sec:B} sont consacrées respectivement 
à l'explication des conditions {\bf B'} et {\bf B}, et la proposition
\ref{prop:vague} est finalement prouvée dans la section \ref{sec:dem}.

\section{La théorie de Breuil}
\label{sec:breuil}

\subsection{Définitions}

\paragraph{Objets d'algèbre linéaire}

On fixe $\kappa$ un entier compris\footnote{Dans la suite de l'exposé, 
on n'utilisera que les valeurs particulières $\kappa = 2$ et $\kappa = 
p-1$.} entre $2$ et $p-1$, ainsi que $E$ une extension finie de $\F_p$
contenue dans $\bar \F_p$. Posons $\tilde S = (k \otimes_{\F_p} E)[u]
/u^{ep}$ (où, rappelons-le, $e = p^r-1$). Pour tout $g \in G$, on note
$\hat g : \tilde S \to \tilde S$ l'unique endomorphisme de $(k 
\otimes_{\F_p} E)$-algèbres qui envoie $u$ sur $(\omega(g) \otimes 1) 
u$. Finalement, soit $\phi : \tilde S \to \tilde S$ défini par 
l'élévation à la puissance $p$ sur $k[u]/u^{ep}$ et l'identité sur $E$.

\begin{deftn}
On note $\EBrModdd^{\kappa-1}$ la catégorie dont les objets sont la 
donnée :
\begin{itemize}
\item[i)] d'un $\tilde S$-module libre $\calM$ ;
\item[ii)] d'un sous-module $\Fil^{\kappa-1} \calM$ de $\calM$ contenant
$u^{e(\kappa-1)} \calM$ ;
\item[iii)] d'un opérateur $\phi$-semi-linéaire $\phi_{\kappa-1} : 
\Fil^{\kappa-1} \calM \to \calM$ dont l'image engendre $\calM$ ;
\item[iv)] d'un opérateur $(k\otimes_{\F_p} E)$-linéaire $N : \calM \to
u \calM$ vérifiant $N(ux) = u N(x) - ux$ pour tout $x \in \calM$
(condition de Leibniz), $u^e N(\Fil^{\kappa-1} \calM) \subset 
\Fil^{\kappa-1} \calM$ et 
$\phi_{\kappa-1}(u^e N(x)) = N(\phi_{\kappa-1}(x))$ pour tout $x \in 
\Fil^{\kappa-1} \calM$ ;
\item[v)] de morphismes $\hat g$-semi-linéaires $[g] : \calM \to \calM$
($g$ parcourant $G$) vérifiant $[\id] = \id$ et $[gh] = [g] \circ [h]$ 
pour tous $g$ et $h$.
\end{itemize}
Les morphismes de cette catégorie sont les applications 
$\tilde S$-linéaires commutant à toutes les structures supplémentaires.
\end{deftn}

\medskip

\noindent
Faisons tout de suite deux remarques importantes. \emph{Primo}, lorsque 
$\kappa = 2$, la donnée du $N$ est automatique : on entend par là que si 
$\calM$ est un $\tilde S$-module libre muni d'un $\Fil^1 \calM$, d'un 
$\phi_1$ et de $[g]$ vérifiant les conditions de la définition, alors il 
existe un et un unique $N$ défini sur ce module qui en fait un objet de 
$\EBrModdd^1$.

\emph{Secundo}, si $\kappa \leq \kappa'$, on peut voir 
$\EBrModdd^{\kappa-1}$ comme une sous-catégorie pleine de 
$\EBrModdd^{\kappa'-1}$ grâce au foncteur pleinement fidèle $\iota : 
\EBrModdd^{\kappa-1} \to \EBrModdd^{\kappa'-1}$ défini comme suit. 
À $\calM \in \EBrModdd^{\kappa-1}$, on associe $\iota(\calM) = \calM$ 
(en tant que $\tilde S$-module) muni de $\Fil^{\kappa'-1} \iota(\calM) = 
u^{e(\kappa'-\kappa)} \Fil^{\kappa-1} \calM$, de $\phi_{\kappa'-1}$ 
défini par $\phi_{\kappa'-1}(u^{e(\kappa'-\kappa)} x) = \phi_{\kappa-1} 
(x)$ et des morphismes $N$ et $[g]$ qui restent inchangés. On remarque
en particulier, que seules les structures qui ont changé de nom ont 
aussi changé de définition ; ceci légitime le léger abus qui consiste à 
noter simplement $\calM$ pour $\iota(\calM)$, ce que nous ferons par
la suite. Signalons finalement que $\iota$ commute au foncteur $\Tst$,
que nous définissons ci-après.

\paragraph{Foncteur vers Galois}

La catégorie $\EBrModdd^{\kappa-1}$ est munie d'un foncteur 
contravariant $\Tst$ vers la catégorie des $E$-représentations
galoisiennes de $G_{K_0}$ notée $\Rep_E(G_{K_0})$. Il est construit à
partir d'un certain anneau de périodes. Contrairement à l'usage, la
définition de ce dernier n'est ici pas vraiment technique, puisqu'il 
s'agit simplement de $\hat A = \O_{\bar K}/p \brac X$ où la notation 
$\brac \cdot$ fait référence à l'algèbre polynômiale à puissances 
divisées. Cet anneau est une $k[u]/u^{ep}$-algèbre grâce au morphisme 
$k[u]/u^{ep} \to \hat A$, $u \mapsto \frac{\pi_1}{1+X}$ où $\pi_1$ est 
une racine $p$-ième de $\pi$ fixée. Il est muni d'un idéal $\Fil^1 \hat 
A$ défini comme celui engendré par une racine $p$-ième de $p$, notée 
$p_1$, et les puissances divisées $\gamma_i(X)$ pour $i \geq 1$. On 
définit un opérateur $\phi_1 : \Fil^1 \hat A \to \hat A$ semi-linéaire 
par rapport au Frobenius sur $\hat A$ par :
$$\phi_1 (p_1) = -1 \quad ; \quad \phi_1(X) = \frac{(1+X)^p-1} p \quad ; 
\quad \phi_1(\gamma_i(X)) = 0, \, i \geq 2$$
où, bien entendu, la fraction est calculée dans $\Z_p[X]$ avant d'être
réduite dans $\hat A$. Pour $1 \leq t \leq p-2$, on note $\Fil^t \hat A$ 
la $t$-ième puissance de $\Fil^1 \hat A$ et on définit un morphisme
additif $\phi_t : \Fil^t \hat A \to \hat A$ en posant $\phi_t(x_1 \cdots 
x_t) = \phi_1(x_1) \cdots \phi_1(x_t)$ pour $x_i$ des éléments de 
$\Fil^1 \hat A$.
En outre, $\hat A$ est muni d'un opérateur de monodromie $N$ qui est 
l'unique morphisme $\O_{\bar K}$-linéaire qui envoie $\gamma_i(X)$ sur 
$(1+X)\gamma_{i-1}(X)$ pour tout $i \geq 1$. Finalement, $\hat A$ est 
muni d'une action de $G_K$ \emph{via} la formule $g(X) = 
\omega(g)(1+X)-1$ valable pour tout $g \in G_{K_0}$. Le foncteur $\Tst$ 
s'obtient alors comme suit :
$$\Tst(\calM) = \hom_{k[u]/u^{ep}, \Fil^{\kappa-1}, \phi_{\kappa-1},N} 
(\calM, \hat A)$$
où la notation signifie que l'on considère les morphismes 
$k[u]/u^{ep}$-linéaires compatibles à $\Fil^{\kappa-1}$, à 
$\phi_{\kappa-1}$ et à $N$. Le groupe $G_{K_0}$ agit sur ce module 
par la formule :
\begin{equation}
\label{eq:galois}
g \cdot f : x \mapsto  g \cdot f([\bar g^{-1}] (x))
\end{equation}
où $g \in G_{K_0}$ et $\bar g$ désigne son image dans $G$. L'action de 
$E$ sur $\Tst(\calM)$ se fait, quant à elle, par l'intermédiaire de son 
action sur $\calM$.

\medskip

Le foncteur ainsi défini $\Tst$ est fidèle. Il est aussi exact dans le 
sens où il transforme les suites exactes courtes dans la catégorie 
$\EBrModdd^{\kappa-1}$ (c'est-à-dire les suites exactes courtes sur les
$\tilde S$-modules sous-jacents qui induisent des suites encore exactes 
sur les $\Fil^{\kappa-1}$) en suites exactes courtes de représentations
galoisiennes. De plus, on a l'égalité de dimension suivante :
$$\dim_E \Tst(\calM) = \dim_{\tilde S} \calM$$
valable pour tout objet $\calM \in \EBrModdd^{\kappa-1}$.

\subsection{Description en rang $1$}

On suppose désormais que $E$ contient un sous-corps isomorphe (non 
canoniquement) à $k$. Il en résulte un isomorphisme d'anneaux $k 
\otimes_{\F_p} E \simeq E^S$ donné par $x \otimes y \mapsto (\tau_i(x) 
y)_{i \in S}$. En considérant les idempotents associés à cette 
décomposition, on déduit que tout module sur $k \otimes_{\F_p} E$ se 
décompose de façon canonique comme une somme directe de $r$ espaces 
vectoriels sur $E$, $k$ agissant sur cette écriture diagonalement 
\emph{via} ses divers plongements dans $E$.

En particulier, si $\calM$ un objet de $\EBrModdd^{\kappa-1}$, on 
peut écrire :
$$\calM = \calM_1 \oplus \calM_2 \oplus \cdots \oplus \calM_r$$
où les $\calM_i$ sont des $E[u]/u^{ep}$-modules libres. De même :
$$\Fil^{\kappa-1} \calM = \Fil^{\kappa-1} \calM_1 \oplus \Fil^{\kappa-1} 
\calM_2 \oplus \cdots \oplus \Fil^{\kappa-1} \calM_r$$
avec $u^{e(\kappa-1)} \calM_i \subset \Fil^{\kappa-1} \calM_i \subset 
\calM_i$. L'action de Frobenius sur $k \otimes_{\F_p} E$ correspond au 
décalage vers la droite sur $E^S$ ; on en déduit que l'opérateur 
$\phi_{\kappa-1}$ induit pour tout $i \in S$ des applications 
$\phi_{\kappa-1} : \Fil^{\kappa-1}
\calM_i \to \calM_{i+1}$ dont l'image engendre tout $\calM_{i+1}$. Du 
fait que les opérateurs $N$ et $[g]$ commutent à l'action de $E$, il 
suit qu'ils stabilisent chacun des $\calM_i$.

\bigskip

Soit maintenant $\calM$ un objet de $\EBrModdd^{\kappa-1}$ qui est de rang $1$ 
comme $\tilde S$-module. Dans la décomposition $\calM = \bigoplus_{i \in 
S} \calM_i$, chacun des $\calM_i$ est libre de rang $1$ sur 
$E[u]/u^{ep}$. Ainsi, pour tout $i$, il existe un unique entier $m_i \in 
[0,e(\kappa-1)]$ tel que $\Fil^{\kappa-1} \calM_i = u^{m_i} \calM_i$.

Notons $e_1$ une base de $\calM_1$. Posons $e_2 = 
\phi_{\kappa-1}(u^{n_1} e_1)$. Comme $\phi_{\kappa-1} : \Fil^{\kappa-1} 
\calM_1 \to \calM_2$ engendre son image, $e_2$ forme nécessairement une 
base de $\calM_2$. Définissons plus généralement par récurrence 
$e_{i+1} = \phi_{\kappa-1}(u^{m_i} e_i)$. Chacun des $e_i$ est une base 
de $\calM_i$ et par conséquent, on a une égalité de la forme $e_{r+1} = a 
e_1$ pour un certain $a$ inversible dans $E[u]/u^{ep}$.

Par ailleurs, si $\lambda$ est un élément inversible de $E[u]/u^{ep}$, 
on constate que modifier $e_1$ en $\lambda e_1$ modifie $e_{r+1}$ en 
$\phi^r(\lambda)\: e_{r+1}$ et donc $a$ en $a \: \frac{\phi^r(\lambda)} 
\lambda$ (où, rappelons-le, $\phi$ agit sur $E[u]/u^{ep}$ en laissant 
fixe $E$ et en envoyant $u$ sur $u^p$). En ajustant correctement 
$\lambda$, on peut donc s'arranger pour avoir $a \in E^\star$.

\medskip

Examinons maintenant l'action de la donnée de descente. Fixons $i \in 
S$. Pour tout $g \in G$, $[g]$ induit un automorphisme de $\calM_i$. 
Ainsi, on peut écrire $[g] e_i = \alpha(g) e_i$ pour une certaine 
fonction $\alpha : \Gal(K/K_0) \to (E[u]/u^{ep})^\times$. La 
compatibilité au produit assure que $\alpha$ est un caractère. Il prend 
ainsi ses valeurs parmi les racines $(p^r-1)$-ièmes de l'unité dans
$(E[u]/u^{ep})^\times$, dont on vérifie facilement qu'elle sont toutes
dans $E^\star$. Ainsi $\alpha$ est un caractère de $G$ à valeurs dans 
$\bar \F_p^\star$ et, en tant que tel, il s'écrit comme une puissance de 
$\omega_i$, disons $\omega_i^{\mu_i}$ où $\mu_i$ est défini modulo 
$(p^r-1)$.

En écrivant à présent la commutation de $[g]$ à $\phi_1$, on obtient la 
relation $\mu_{i+1} \equiv p (\mu_i + m_i) \pmod{p^r-1}$. Une 
combinaison linéaire judicieuse de ces relations permet d'éliminer les
$\mu_i$ et conduit à la congruence :
\begin{equation}
\label{eq:conmi}
p^r m_i + p^{r-1} m_{i+1} + \cdots + p^2 m_{i+r-2} + p m_{i+r-1} \equiv 
0 \pmod {p^r-1}.
\end{equation}

Déterminons finalement l'opérateur de monodromie. La relation de Leibniz
montre que $u^{m_i}$ divise $N(u^{m_i} e_i)$. Comme ce dernier est 
élément de $\calM_i$, il appartient nécessairement à $\Fil^{\kappa-1} 
\calM_i$. Ainsi $\phi_{\kappa-1} (u^e N(u^{m_i} e_i)) = u^{ep} 
\phi_{\kappa-1} (N(u^{m_i} e_i)) = 0$, et donc $N(e_{i+1}) = 0$. En 
résumé, on vient de prouver la proposition suivante :

\begin{prop}
\label{prop:rang1}
Soit $\calM$ un objet de $\EBrModdd$ de rang $1$. Alors il existe des 
éléments $e_i \in \calM_i$, des entiers $m_i$ compris entre $0$ et 
$e(\kappa-1)$, des entiers $\mu_i$ définis modulo $(p^r-1)$ et un 
élément $a \in E^\star$ tels que :
\begin{itemize}
\item[i)] l'élément $e_i$ forme une base de $\calM_i$ ;
\item[ii)] l'élément $u^{m_i} e_i$ engendre $\Fil^{\kappa-1} \calM_i$ ;
\item[iii)] $\phi_{\kappa-1}(u^{m_i} e_i) = e_{i+1}$ pour $1 \leq i \leq 
r-1$ et $\phi_{\kappa-1}(u^{m_r} e_r) = a e_1$ ;
\item[iv)] $\mu_{i+1} \equiv p (\mu_i + m_i) \pmod{p^r-1}$ ;
\item[v)] pour tout $g \in G$, $[g] (e_i) = \omega_i^{\mu_i}(g) e_i$ ;
\item[vi)] $N(e_i) = 0$.
\end{itemize}
Dans ces conditions, les entiers $m_i$ vérifient automatiquement la
congruence (\ref{eq:conmi}) et nous notons :
$$\mu_{\f,i} = \frac{p^r m_i + p^{r-1} m_{i+1} + \cdots + p^2 
m_{i+r-2} + p m_{i+r-1}}{p^r-1}.$$
Réciproquement, l'objet défini précédemment est dans la catégorie
$\EBrModdd$.
\end{prop}

\medskip

\noindent
{\it Remarque.} Dans la proposition précédente, les entiers $m_i$, les
classes de congruence des $\mu_i$ et l'élément $a$ sont uniquement 
determinés. Ce n'est par contre pas le cas des $e_i$ qui peuvent être 
changés en $e'_i = \lambda e_i$ pour n'importe quel $\lambda \in 
E^\star$ (le même pour chaque $e_i$). On peut vérifier de surcroît qu'il 
s'agit là du seul degré de liberté autorisé.

\begin{prop}
\label{prop:gal1}
Avec les notations de la proposition \ref{prop:rang1}, l'action de 
$G_{K_0}$ sur $\Tst(\calM)$ se fait par l'intermédiaire du caractère :
$$\lambda_a \cdot \omega_i^{(\kappa-1)(1 + p + p^2 + \cdots + p^{r-1}) - 
(\mu_i + \mu_{\f,i})}$$
pour tout $i\in S$. Ici, $\lambda_a$ est le caractère non ramifié 
qui envoie un Frobenius géométique sur $a$.
\end{prop}

\noindent
{\it Remarque.} En particulier, $\omega_i^{\mu_i + \mu_{\f,i}}$ ne 
dépend pas du choix de $i$, ce qui signifie que l'on a la congruence $p 
(\mu_i + \mu_{\f,i}) \equiv \mu_{i+1} + \mu_{\f,i+1} \pmod{p^r-1}$ pour 
tout $i$. Ceci peut se vérifier aisément par ailleurs.

\begin{proof}
Nous donnons simplement quelques idées de la preuve. On utilise bien 
entendu la formule $\Tst(\calM) = \hom_{k[u]/u^{ep}, \Fil^{\kappa-1}, 
\phi_{\kappa-1}, N} (\calM, \hat A)$.

Notons $n$ le degré de $E$ sur $k$ et supposons pour simplifier $a$ 
engendre cette extension. Alors la famille des $a^j e_i$ ($0 \leq j
\leq n-1$, $i \in S$) est une $(k[u]/u^{ep})$-base de $\calM$ et un
élément de $\Tst(\calM)$ est entièrement déterminé par les images de ces
vecteurs, images qui sont soumises à certaines relations. La première
de celle-ci est $N(x_{ij}) = 0$ (puisque $f$ doit être compatible à 
l'action de $N$), c'est-à-dire $x_{ij} \in \O_{\bar K}/p$. Les autres
relations s'obtiennent en écrivant la compatibilité à $\Fil^{\kappa-1}$ 
et $\phi_{\kappa-1}$. 

Par ailleurs, on sait que $\Tst(\calM)$ est un $E$-espace vectoriel
de dimension $1$ et donc qu'il a même cardinal que $E$, c'est-à-dire
$p^{nr}$. Il suffit donc, pour déterminer $\Tst(\calM)$, de trouver 
$p^{nr}$ solutions distinctes pour les $x_{ij}$ qui s'intuitent en fait
assez facilement une fois que les équations sur les $x_{ij}$ ont été
écrites.

Il ne reste alors plus qu'à comprendre l'action de Galois.
\end{proof}

\section{La condition B' : schémas en groupes}
\label{sec:Bp}

\subsection{L'équivalence de Breuil}

Appelons \emph{$E$-groupe sur $\O_K$} un schéma en $E$-vectoriels fini 
et plat sur $\O_K$. Une donnée de descente (de $K$ à $K_0$) sur un 
$E$-groupe $\G$ sur $\O_K$ est la donnée, pour tout élément $g \in 
G$, d'un morphisme $[g] : \G \to \G$ rendant le diagramme suivant 
commutatif :
$$\xymatrix @R=13pt @C=50pt {
\G \ar[d] \ar[r]^-{[g]} & \G \ar[d] \\
\Spec(\O_K) \ar[r]_-{\Spec(g)} & \Spec(\O_K) }$$
et tel que le morphisme déduit $\G \to {}^g \G$ (où ${}^g \G$ est 
déduit de $\G$ par le changement de base $\Spec(g)$) soit un
morphisme de $E$-groupes, le tout étant soumis aux relations $[\id] =
\id$ et $[gh] = [g] \circ [h]$ pour $g$ et $h$ dans $G$.

Notons $\EGrdd$ la catégorie des $E$-groupes sur $\O_K$ munis d'une 
donnée de descente. Le théorème fondamental suivant (voir \cite{breuil}, 
\cite{bcdt} et éventuellement \cite{savitt}) permet de comprendre 
comment les catégories introduites dans la section \ref{sec:breuil} 
interviennent dans notre propos.

\begin{theo}[Breuil]
Supposons $p > 2$.
Il existe une anti-équivalence\footnote{Dans \cite{gee}, Gee travaille
plutôt avec le foncteur covariant $\Gr$ défini par $\Gr(\calM) = 
\Gr^\star(\calM)^\vee$ où \og $\vee$ \fg\ désigne le dual de Cartier.
Cependant, dans (la preuve de) la proposition 3.3.1, il est 
contraint d'utiliser la version contravariante, et c'est pourquoi nous 
présentons celle-ci dans cette note.} de catégories $\Gr^\star : 
\EBrModdd \to \EGrdd$ qui rend le diagramme suivant commutatif :
$$\xymatrix @R=10pt @C=50pt {
\EBrModdd^1 \ar[dd]_-{\Gr^\star}^-{\sim} \ar[rd]^-{\Tst} \\
& \Rep_E(G_{K_0}) \\
\EGrdd \ar[ru] }$$
où le foncteur $\EGrdd \to \Rep_E(G_{K_0})$ est celui donné par les $\bar 
K$-points.
\end{theo}

\noindent
{\it Remarque.} Dans le théorème précédent, la catégorie de modules
considérée correspond à $\kappa = 2$. En particulier, on est dans
la situation où l'opérateur $N$ est défini de façon automatique.

\subsection{Schémas en groupes de classe $J$}

Soit $\G$ un $E$-groupe sur $\O_K$ avec donnée de descente et $\calM$ 
son module associé. On suppose que $\calM$ est de rang $1$, et donc 
qu'il a la forme donnée par la proposition \ref{prop:rang1}. Appelons 
$m_i$, $\mu$ et $a$ les invariants numériques qui interviennent. La 
donnée des $m_i$ permet de caractériser les groupes étales et ceux de 
type multiplication : exactement $\G$ est étale (resp. de type 
multiplicatif) si, et seulement si tous les $m_i$ sont égaux à 
$e$ (resp. sont nuls). Pour $J$ un sous-ensemble de $S$, nous 
introduisons suivant Gee, une notion intermédiaire entre ces deux 
extrêmes.

\begin{deftn}
\label{def:classeJ}
Soit $J \subset S$. Avec les notations précédentes, on dit que $\G$ 
est de \emph{classe $J$} si\footnote{La formule que nous donnons diffère 
en deux endroits de celle que l'on peut lire dans \cite{gee}. Tout 
d'abord, nous écrivons $\bar \1_J$ à la place $\1_J$, mais cela est 
simplement dû au fait que nous utilisons le foncteur $\Gr^\star$ (et
pas $\Gr$). Le second écart est que nous appliquons cette fonction non 
pas à $i$, mais à $i+1$ : cela fait par contre une différence 
considérable et c'est ici que réside l'erreur dans \cite{gee} que nous 
mentionnions au début de cette note.} $m_i = e \: \bar 
\1_J(i+1)$.
\end{deftn}

\noindent
{\it Remarque.} Comme $e = p^r-1$, les $m_i$ donnés par la définition 
vérifient toujours la congruence (\ref{eq:conmi}).

\begin{prop}
\label{prop:grpJ}
Fixons $J$ un sous-ensemble de $S$ ainsi qu'un caractère $\psi : 
G_{K_0} \to E^\star$ trivial sur $I_K$. Alors, il existe un unique 
$E$-groupe de classe $J$ sur $\O_K$ (avec donnée de descente) dont la 
représentation galoisienne associée correspond à $\psi$.
\end{prop}

\begin{proof}
Comme $\psi$ est trivial sur $I_K$, il s'écrit sous la forme $\lambda_b 
\omega_0^n$ où $\lambda_b$ est le caractère non ramifié qui envoie un 
Frobenius géométique sur $b$.
Le résultat découle alors des propositions \ref{prop:rang1} et 
\ref{prop:gal1} : avec leurs notations, il suffit de choisir :
$$\mu_i \equiv - p^i n + \sum_{j=1}^r p^{r-j} \1_J(i+j+1) \pmod{p^r-1}$$
et $a=b$.
\end{proof}

\noindent
{\it Remarque.} D'après la théorie du corps de classe, l'hypothèse de
trivialité sur $I_K$ est en fait automatique.

\subsection{Explication de la condition B'}

Prenons pour $E$ un corps suffisamment grand pour que $\rho$ se 
factorise par $\GL_2(E)$. On suppose qu'il existe un sous-ensemble $J$ 
de $S$ pour lequel la condition {\bf A} est satisfaite avec $a_i = 0$ et
$2 \leq b_i \leq p-2$ pour tout $i$. D'après la proposition 
\ref{prop:grpJ}, il existe un unique $E$-groupe $\G'$ (resp. $\G''$) 
avec donnée de descente de classe $J$ (resp. de classe $S \backslash J$) 
dont la représentation galoisienne associée s'identifie à $E(\psi')$ 
(resp. $E(\psi'')$) \emph{via} un isomorphisme $f'$ (resp. $f''$). La 
condition {\bf B'} s'exprime alors comme suit :

\medskip

\begin{itemize}
\item[\it Condition {\bf B'}] : Il existe un $E$-groupe $\G$ avec donnée 
de descente qui s'insère dans une suite exacte courte $0 \to \G' \to \G 
\to \G'' \to 0$ et un isomorphisme $G_{K_0}$-équivariant $f : \G(\bar K) 
\to E^2(\rho)$ (où $E^2(\rho)$ est la représentation de dimension $2$ 
donnée par $\rho$) s'insérant dans un diagramme commutatif :
$$\xymatrix @C=40pt @R=13pt {
0 \ar[r] & \G'(\bar K) \ar[r] \ar[d]_-{f'} & \G(\bar K) \ar[r] \ar[d]_-{f} 
& \G''(\bar K) \ar[r] \ar[d]^-{f''} & 0 \\
0 \ar[r] & E(\psi') \ar[r] & E^2(\rho) \ar[r] & E(\psi'') \ar[r] & 0
}$$
\end{itemize}

\section{La condition B : relevés cristallins}
\label{sec:B}

\subsection{Représentations cristallines}

Posons $\O_L = W(E)$ et $L = \O_L[1/p]$. Considérons le produit 
tensoriel $K_0 \otimes_{\Q_p} L$ et munissons-le d'un endomorphisme 
$\phi$ agissant comme le Frobenius sur $K_0$ et comme l'identité sur 
$L$. Comme précédemment, on a un isomorphisme d'anneaux $K_0 
\otimes_{\Q_p} L \simeq L^S$ et l'opérateur $\phi$ correspond par cet 
isomorphisme au décalage vers la droite.

\begin{deftn}
Un \emph{$L$-module filtré} est un $(K_0 \otimes_{\Q_p} L)$-module libre 
$D$ muni d'un opérateur $\phi$-semi-linéaire $\phi : D \to D$ et d'une 
filtration par des $(K_0 \otimes_{\Q_p} L)$-sous-modules (pas 
nécessairement libres) $\Fil^t D$ exhaustive et séparée.
\end{deftn}

Soit $D$ un $L$-module filtré. Il lui est associé deux invariants 
numériques. Le premier est son nombre de Newton, noté $t_N(D)$, et 
est défini comme la pente de $\phi$ sur la puissance extérieure maximale
de $D$ considéré comme $K_0$-espace vectoriel. Le second est son nombre
de Hodge, noté $t_H(D)$, et est défini comme l'unique saut de la filtration 
sur la même puissance extérieure maximale. On dit que $D$ est 
\emph{faiblement admissible} si $t_H(D) = t_N(D)$ et si pour tout $D' 
\subset D$ stable par $\phi$ et muni de la filtration \og intersection 
\fg, on a $t_H(D') \leq t_N(D')$. Fontaine construit un foncteur
pleinement fidèle :
$$\Vcris : \{L\text{-module filtrés faiblement admissibles}\} \to 
\Rep_L(G_{K_0})$$
défini par $\Vcris(D) = \hom_{K_0, \Fil^\cdot, \phi} (D, B_\cris)$ 
où $B_\cris$ est un anneau de périodes dont la définition est par
exemple donnée dans \cite{fontaine}. Le foncteur $\Vcris$ est 
de plus exact dans le sens où il transforme suites exactes courtes en 
suites exactes courtes, une suite de $L$-modules filtrés étant dite
exacte si elle induit pour tout $t$ une suite exacte sur les $\Fil^t$.
Finalement, l'image essentielle de $\Vcris$ constitue par 
définition ce que l'on appelle les \emph{représentations cristallines}.

\subsection{Description en rang $1$}

De l'isomorphisme $K_0 \otimes_{\Q_p} L \simeq L^S$, il 
résulte que tout module sur cet anneau s'écrit canoniquement comme une 
somme directe de $r$ espaces vectoriels sur $L$. En particulier, si $D$ 
est un $L$-module filtré, on peut écrire :
\begin{equation}
\label{eq:decomp}
D = D_1 \oplus D_2 \oplus \cdots \oplus D_r
\quad \text{ et } \quad 
\Fil^t D = \Fil^t D_1 \oplus \Fil^t D_2 \oplus \cdots \oplus 
\Fil^t D_r
\end{equation}
pour tout $t$. Les $\Fil^t D_i$ forment alors une filtration 
décroissante exhaustive et séparée de $D_i$. De plus le Frobenius induit 
par restriction des applications $\phi : D_i \to D_{i+1}$.

\medskip

Supposons à présent que $D$ est de rang $1$. Les $D_i$ qui interviennent
dans la décomposition précédente sont des $L$-espaces vectoriels de
dimension $1$. Soient $e'_1$ une base de $D_1$, $e'_2 = \phi(e'_1)$,
$e'_3 = \phi(e'_2)$, \emph{etc}. Comme $e'_1$ et $e'_{r+1}$ sont tous
les deux des bases de $D_1$, on a $e'_{r+1} = \lambda e'_1$ pour un
certain $\lambda \in L^\star$. Par ailleurs comme $D_i$ étant de
dimension $1$, il existe un unique entier $m_i$ tel que $\Fil^t D_i =
D_i$ si $t \leq m_i$ et $\Fil^t D_i = 0$ sinon. La condition de faible
admissibilité se traduit ici simplement par l'égalité :
$$v(\lambda) = m_1 + m_2 + \cdots + m_r$$
où $v$ est la valuation sur $L$ normalisée par $v(p) = 1$. Ainsi, on 
peut écrire $\lambda = p^{m_1 + \cdots + m_r} x$ où $x$ est un élément
inversible de $\O_L$. On a ainsi obtenu une description complète de $D$. 
On peut la rendre un peu plus symétrique comme le résume la proposition 
suivante :

\begin{prop}
\label{prop:ratrang1}
Soit $D$ un $L$-module filtré faiblement admissible de rang $1$ sur $K_0 
\otimes_{\Q_p} L$. Alors, il existe des entiers $m_i$, un élément $x$ 
inversible dans $\O_L$ et des $e_i \in D$ tels que :
\begin{itemize}
\item[i)] pour tout $i \in S$, $e_i$ forme une base de $D_i$ ;
\item[ii)] pour tout $i \in S$, $\Fil^t D_i = D_i$ si $t \leq m_i$ et 
$\Fil^t D_i = 0$ sinon ;
\item[iii)] pour $1 \leq i \leq r-1$, $\phi(e_i) = p^{m_i} e_{i+1}$ et
$\phi(e_r) = x p^{m_r} e_1$.
\end{itemize}
\end{prop}

\begin{proof}
Il suffit de poser $e_i = p^{m_i + m_{i+1} + \cdots + m_r} e'_i$.
\end{proof}

\noindent
{\it Remarque.} Les entiers $m_i$ et l'élément $a$ sont uniquement 
déterminés. De plus, les $m_i$ sont reliés à des invariants plus 
usuels, puisque ce sont les opposés des poids de Hodge-Tate de la 
$\Q_p$-représentation associée à $D$ par le foncteur $\Vcris$.

\subsection{Théorie de Fontaine-Laffaille}
\label{subsec:fl}

L'un des buts de la théorie de Fontaine-Laffaille est de décrire les
réseaux stables par l'action de Galois dans les représentations
cristallines en terme de modules filtrés. Pour y parvenir, la méthode
consiste à définir des équivalents de ces réseaux dans les modules
filtrés faiblement admissibles discutés précédemment. On aura cependant
besoin d'une hypothèse supplémentaire sur la filtration : on demande 
qu'elle soit concentrée en bas degré. Plus précisément, Fontaine et
Laffaille posent la définition suivante :

\begin{deftn}
\label{def:reseaufl}
Soit $D$ un $L$-module filtré. On suppose $\Fil^0 D = D$ et $\Fil^{p-1} 
D = 0$. Un \emph{$\O_L$-réseau} de $D$ est la donnée d'un 
sous-$(\O_{K_0} \otimes_{\Z_p} \O_L)$-module $M$ de $D$ vérifiant :
\begin{itemize}
\item[i)] $M[1/p] = D$ ;
\item[ii)] pour tout $t$, l'application $\phi_t = \frac{\phi}{p^t}$
envoie $\Fil^t M = M \cap \Fil^t D$ dans $M$ ;
\item[iii)] $\sum_{t=0}^{p-2} \phi_t (\Fil^t M) = M$.
\end{itemize}
\end{deftn}

\noindent
{\it Remarques.} On n'a fait aucun hypothèse de faible admissibilité 
dans la définition précédente. En réalité, l'existence d'un réseau 
possédant ces propriétés est équivalente à la faible admissibilité de 
$D$.

Le prototype de modules filtrés de rang $1$ donné par la proposition
\ref{prop:ratrang1} admet le réseau évident $M = \sum_{i=1}^r \O_L e_i$.

\bigskip

Par ailleurs, $B_\cris$ contient lui aussi une structure entière qui 
consiste en un sous-anneau $A_\cris$ (voir \cite{fontaine} pour une
définition) qui permet de relier les réseaux définis précédemment aux 
réseaux galoisiens. Précisément, considérons $D$ un $L$-module filtré
faiblement admissible tel que $\Fil^0 D = D$ et $\Fil^{p-1} D = 0$ et
notons $V = V_\cris^\star(D)$. On associe alors à tout $\O_L$-réseau 
$M \subset D$ la représentation :
$$\Tcris(M) = \hom_{\O_{K_0}, \Fil^\cdot, \phi} (M, A_\cris)$$
qui est un réseau (au sens usuel) dans $V$. Fontaine et 
Laffaille montre que cette recette établit une bijection entre les 
réseaux de $D$ et les réseaux de $V$ stables par Galois. De plus, ils 
prouvent que le foncteur $\Tcris$ est lui aussi exact.

\subsection{Réduction modulo $p$}
\label{subsec:modp}

Considérons $D$ un $L$-module filtré. À un $\O_L$-réseau $M \subset D$, 
il est possible d'associer un objet de $\EBrModdd^{p-2}$ qui calcule
la réduction modulo $p$ de $\Tcris(M)$. Ceci se fait par les formules 
suivantes :
$$\begin{array}{c}
\calM = \tilde S \otimes_{\O_{K_0} \otimes_{\Z_p} \O_L} M 
\quad ; \quad \Fil^{p-2} \calM = \displaystyle \sum_{t=0}^{p-2} 
u^{et} \tilde S \otimes_{\O_{K_0} \otimes_{\Z_p} \O_L} \Fil^{p-2-t} 
M \\
\phi_{p-2} = \sum \phi_t \otimes \phi_{p-2-t} \quad ; \quad N = N 
\otimes \id \quad ; \quad [g] = \hat g \otimes \id \text{ pour tout } g 
\in G
\end{array}$$
où par définition $\phi_t : u^{et} \tilde S \to \tilde S$ est l'unique
application $\phi$-semi-linéaire qui envoie $u^{et}$ sur $1$. Le fait
que $\calM$ calcule la réduction modulo $p$ de $\Tcris(M)$ signifie
que l'on a une identification canonique $\Tcris(M)/p\Tcris(M) \simeq 
\Tst(\calM)$.

En particulier, considérons $\chi : G_{K_0} \to L^\star$ le caractère 
qui correspond \emph{via} $\Vcris$ à l'objet $D$ de la proposition 
\ref{prop:ratrang1}. Comme $G_{K_0}$ est un groupe compact, il prend ses 
valeurs dans $\O_L^\times$, et on peut considérer $\psi$ sa réduction 
modulo $p$ ; c'est un caractère à valeurs dans $E^\star$ et un calcul 
direct à partir de ce qui précède montre qu'il s'exprime comme suit :
\begin{equation}
\label{eq:galres}
\psi = \lambda_{\bar x} \cdot \omega_1^{m_1} \omega_2^{m_2} \cdots 
\omega_r^{m_r}
\end{equation}
où $\lambda_{\bar x}$ est le caractère non ramifié de $G_{K_0}$ qui 
envoie un Frobenius géométrique sur l'image de $x \in \O_L^\times$ dans 
$E^\star$.

\subsection{Explication de la condition B}

On considère $E$ suffisamment grand pour que $\rho$ se factorise par 
$\GL_2(E)$. On suppose qu'il existe un sous-ensemble $J$ de $S$ pour 
lequel la condition {\bf A} est satisfaite avec $a_i = 0$ et $2 \leq b_i 
\leq p-2$ pour tout $i \in S$. D'après la formule (\ref{eq:galres}), il 
existe un unique $L$-module filtré $D'$ (resp. $D''$) de rang $1$ tel 
que :
\begin{itemize}
\item[$\bullet$] la réduction modulo $p$ d'un réseau galoisien à 
l'intérieur de $\Vcris(D)$ est isomorphe à $E(\psi')$ (resp. 
$E(\psi'')$) par un morphisme noté $f'$ (resp $f''$) ;
\item[$\bullet$] les invariants associés à $D'$ (resp. $D''$) par la 
proposition \ref{prop:ratrang1} sont tels $m_i = b_i \1_J(i)$ (resp. 
$m_i = b_i \bar \1_J(i)$) et $x$ est un représentant de Teichmüller.
\end{itemize}

\bigskip

\noindent
La condition {\bf B} s'exprime alors comme suit :

\medskip

\begin{itemize}
\item[\it Condition {\bf B}] : Il existe un $\O_L$-réseau $M$ (inclus 
dans un $L$-module filtré $D$) qui s'insère dans une suite exacte courte 
$0 \to M'' \to M \to M \to 0$ et un isomorphisme $G_{K_0}$-équivariant 
$f : \Tcris(M) \to E^2(\rho)$ (où $E^2(\rho)$ est la représentation de 
dimension $2$ donnée par $\rho$) s'insérant dans un diagramme 
commutatif :
$$\xymatrix @C=40pt @R=13pt {
0 \ar[r] & \Tcris(M') \ar[r] \ar[d]_-{f'} & \Tcris(M) \ar[r] 
\ar[d]_-{f} & \Tcris(M'') \ar[r] \ar[d]^-{f''} & 0 \\
0 \ar[r] & E(\psi') \ar[r] & E^2(\rho) \ar[r] & E(\psi'') \ar[r] & 0
}$$
\end{itemize}

\bigskip

\noindent
{\it Remarque.} Il est évidemment possible d'utiliser l'équivalence de
catégories $\Vcris$ (ou $\Tcris$) pour exprimer la condition {\bf 
B} simplement en termes de représentations cristallines.

\section{Preuve de la proposition 3.3.1 de \cite{gee}}
\label{sec:dem}

On conserve l'extension $E/\F_p$ suffisamment grande pour que $\rho$ se
factorise par $\GL_2(E)$. On fixe des entiers $b_i$ compris entre $2$ et 
$p-2$ (en particulier, on notera que cela implique $p \geq 5$), et un 
sous-ensemble $J$ de $S$. On suppose que $\rho$ vérifie la 
condition {\bf A} avec $a_i = 0$ pour tout $i \in S$ (et les $b_i$ 
précédents). Le but est donc de 
prouver que les conditions {\bf B} et {\bf B'} expliquées dans les 
sections précédentes sont équivalentes. Notons pour cela $M'$ et 
$M''$ les 
$\O_L$-réseaux qui apparaissent dans la définition de la condition {\bf 
B'}. Notons également $\G'$ et $\G''$ les $E$-groupes avec donnée de 
descente qui apparaissent dans la définition de la condition {\bf B} 
et $\calM'$ et $\calM''$ les objets de $\EBrModdd^1$ associés.

\medskip

Si $x$ est un élément de $E$, désignons par $\tilde x \in \O_L$ son 
représentant de Teichmüller, et si $x$ est un élément, appelons 
$\underline x$ la fonction de domaine $S$ qui associe $x$ à $r$ et $1$ 
aux autres éléments. Avec ces notations, il existe $a$ et $b$ dans $E$ 
tels que l'on ait les descriptions explicites qui suivent :
$$\left\{\begin{array}{l}
\calM'_i = E[u]/u^{ep} \cdot \bar f_i \\
\Fil^1 \calM'_i = u^{e \bar \1_J(i+1)} \calM'_i \\
\phi_1(u^{e \bar \1_J(i+1)} \bar f_i) = \underline a (i) \bar f_{i+1} 
\\{}
[g] (\bar f_i) = \omega_i^{\mu'_i}(g) \bar f_i
\end{array}\right.
\qquad 
\left\{\begin{array}{l}
\calM''_i = E[u]/u^{ep} \cdot e_i \\
\Fil^1 \calM''_i = u^{e \1_J(i+1)} \calM''_i \\
\phi_1(u^{e \1_J(i+1)} e_i) = \underline b (i) e_{i+1} \\{}
[g] (e_i) = \omega_i^{\mu''_i}(g) e_i
\end{array}\right.$$
où :
$$\mu'_i = \sum_{j=1}^r p^{r-j} \big[\1_J (i+j+1) - b_{i+j} 
\1_J(i+j)\big] \qquad \mu''_i = \sum_{j=1}^r p^{r-j} \big[\bar \1_J 
(i+j+1) - b_{i+j} \bar \1_J(i+j)\big]$$
et :
$$\left\{ \begin{array}{l}
M'_i = \O_L \cdot \bar F_i \\
\Fil^t D'_i = D'_i \text{ pour } t \leq b_i \1_J(i) \\
\Fil^t D'_i = 0 \text{ pour } t > b_i \1_J(i) \\
\phi(\bar F_i) = \underline {\tilde a}(i) p^{b_i \1_J(i)} \bar F_{i+1} 
\end{array} \right.
\quad 
\left\{ \begin{array}{l}
M''_i = \O_L \cdot E_i \\
\Fil^t D''_i = D''_i \text{ pour } t \leq b_i \bar \1_J(i) \\
\Fil^t D''_i = 0 \text{ pour } t > b_i \bar \1_J(i) \\
\phi(E_i) = \underline {\tilde b}(i) p^{b_i \bar \1_J(i)} E_{i+1}
\end{array} \right.$$

\subsection{Modules associés aux relevés cristallins}
\label{subsec:demrel}

La proposition suivante donne la forme générale d'un $M$ qui intervient 
dans la condition {\bf B} :

\begin{prop}
\label{prop:descM}
Il existe, pour tout $i$, une $\O_L$-base $(E_i, F_i)$ de $M_i$ et 
des $\Lambda_i \in \O_L$ avec $\Lambda_i = 0$ si $i \not\in J$ et
\begin{itemize}
\item[i)] Pour tout $i \in S$, $\Fil^t M_i = M_i$ pour $t \leq 0$ et
$\Fil^t M_i = 0$ pour $t > b_i$.

\item[ii)] Si $i \in J$, $\Fil^t M_i = \O_L F_i$ pour $0 < t 
\leq b_i$.

Si $i \not\in J$, $\Fil^t M_i = \O_L E_i$ pour $0 < t \leq b_i$.

\item[iii)] $\phi(E_i) = \underline {\tilde b}(i) p^{b_i \bar 
\1_J(i)} E_{i+1}$ et $\phi(F_i) = \underline {\tilde a} (i) p^{b_i
\1_J(i)} (F_{i+1} - \Lambda_{i+1} E_{i+1})$.
\end{itemize}

De plus, le morphisme $N'' \to N$ (resp. $N \to N'$) est donné sur cette 
description par $E_i \mapsto E_i$ (resp. $E_i \mapsto 0$, $F_i \mapsto 
\bar F_i$).
\end{prop}

\begin{proof}
Bien sûr, les $E_i$ sont ceux qui proviennent de $M''$. La condition sur 
la filtration détermine les $F_i$ lorsque $i \in J$. Par ailleurs, la 
condition sur $\phi$, couplée au fait que $\Lambda_i$ s'annule pour $i 
\not\in J$, détermine la valeur de $F_i$ en fonction de $F_{i-1}$ 
lorsque $i \not\in J$. Ainsi lorsque $J \neq \emptyset$, tous les $F_i$ 
sont déterminés et il est immédiat de vérifier qu'ils satisfont toutes
les propriétés de la proposition.

Supposons $J$ vide. Soit $F'_1$ un relevé quelconque de $\bar F_1$. 
Définissons $F'_2 = \phi(F'_1)$, $F'_3 = \phi(F'_2)$ et ainsi de suite. 
Hélas, on n'a pas nécessairement comme cela $F'_{r+1} = \tilde a F'_1$.
Toutefois, la réduction de cette égalité dans $N'$ est vérifiée et donc 
on a $F'_{r+1} = \tilde a F'_1 + \beta E_1$ pour un certain $\beta \in 
\O_L$. Les éléments :
$$F_1 = F'_1 + \frac{\beta}{\tilde a - \tilde b \: p^{b_1 + \cdots + 
b_r}} E_1,$$
$F_2 = \phi(F_1)$, \emph{etc}. répondent alors à la question (notez que 
le dénominateur est bien inversible car $b_1 + \cdots + b_r \geq 2r \geq 
1$).
\end{proof}

\subsection{Modules associés aux $E$-groupes $\G$}
\label{subsec:demgrp}

On souhaite maintenant parvenir à une description explicite analogue 
pour les objets qui interviennent dans la condition {\bf B'}. La donnée 
d'un $E$-groupe (avec donnée de descente) qui s'insère dans une suite
exacte $0 \to \G' \to \G \to \G'' \to 0$ est équivalente à la donnée
d'un objet $\calM \in \EBrModdd^1$ s'insérant dans la suite exacte 
$0 \to \calM'' \to \calM \to \calM' \to 0$. Cherchons à déterminer
la forme de ces $\calM$.

\paragraph{Choix d'un relevé compatible à la donnée de descente}

Fixons $g$ un générateur de $G$. Choisissons $f_1 \in \calM$ un relevé 
de $\bar f_1$ appartenant à l'image de $\phi_1$. On a une écriture de la 
forme :
$$[g] (f_1) = \omega_1^{\mu'_1}(g) f_1 + s e_1$$
pour un certain $s \in E[u^p]/u^{ep}$. Remplacer $f_1$ par $f_1 + 
\delta u^{pn} e_1$ modifie $s$ en $s + \delta [\omega_1^{\mu''_1+pn}(g) 
- \omega_1^{\mu'_1}(g)] u^{pn}$. Ainsi dès que $pn \not\equiv \mu'_1 - 
\mu''_1 \pmod{p^r-1}$, on peut choisir $\alpha$ de sorte à éliminer le 
terme en $u^{pn}$ dans $s$. Par conséquent, si l'on désigne par $n$ 
l'unique entier de $[0,p^r-1[$ congru à $\frac{\mu'_1 - \mu''_1} p$ 
modulo $(p^r-1)$, on peut choisir $f_1$ de sorte à avoir :
\begin{equation}
\label{eq:gfi}
[g] (f_1) = \omega_1^{\mu'_1}(g) f_1 + \alpha(g) u^{pn} e_1
\end{equation}
avec $\alpha(g) \in E$. Libérons $g$. Utilisant la relation 
$[gh] = [g] \circ [h]$, on s'assure que l'égalité (\ref{eq:gfi}) est 
valable pour tout $g \in G$ pour une certain fonction $\alpha : G \to E$ 
soumise à la relation :
$$\alpha(gh) = \omega_1^{\mu'_1}(g) \alpha(h) + \omega_1^{\mu'_1}(h) 
\alpha(g)$$
pour tous $g$ et $h$ dans $G$. Autrement dit la fonction $\beta = \alpha 
\: \omega_i^{-\mu'_1}$ vérifie $\beta(gh) = \beta(g) + \beta(h)$. Comme 
$G$ est d'exposant $p^r-1$ premier à $p$ et que $E$ est de 
caractéristique $p$, la seule solution est d'avoir $\beta = 0$, et par 
suite $\alpha = 0$.

Au final, on a prouvé le lemme suivant :

\begin{lemme}
\label{lem:galois}
Il existe $f_1 \in \phi_1(\Fil^1 \calM)$ relevant $\bar f_1$ et 
vérifiant :
$$[g](f_1) = \omega_1^{\mu'_1} f_1$$
pour tout $g \in G$. De plus, $f_1$ est défini à addition
près d'un élément de la forme $\alpha u^{pn} e_1$ pour $\alpha \in E$
et $n$ l'unique entier de $[0, p^r-1[$ congru à $\frac{\mu'_1 - \mu''_1} 
p$ modulo $(p^r-1)$.
\end{lemme}

\paragraph{Forme de la filtration}

Rappelons que la suite $0 \to \Fil^1 \calM' \to \Fil^1 \calM \to \Fil^1 
\calM'' \to 0$ est exacte. Ainsi $\Fil^1 \calM_1$ est engendré par 
$u^{e \1_J(2)} e_1$ et par un relevé de $u^{e \bar \1_J(2)} \bar f_1$ qui 
appartient à $\Fil^1 \calM_1$.

Supposons dans un premier temps que $2 \not\in J$. Un relevé de $u^e 
\bar f_1$ est alors par exemple $u^e f_1$ : c'est bien un élément de 
$\Fil^1 \calM_1$ puisque $u^e \calM \subset \Fil^1 \calM$. Ainsi $\Fil^1 
\calM_1$ est engendré par $e_1$ et $u^e f_1$.

Supposons maintenant que $2 \in J$. Un relevé dans $\Fil^1 \calM_1$ de 
$\bar f_1$ s'écrit $f_1 + s e_1$ avec $s \in E[u]/u^{ep}$. Ainsi 
$\Fil^1 \calM_1$ est engendré par $f_1 + s e_1$ et $u^e e_1$. On peut 
évidemment choisir $s$ de degré strictement inférieur à $e$. Calculons :
$$[g](f_1 + s e_1) = \omega_1^{\mu'_1}(g) f_1 + \hat g(s) 
\omega_1^{\mu''_1}(g) e_1 = \omega_1^{\mu'_1}(g) (f_1 + s e_1) + \cro{\hat 
g(s) \omega_1^{\mu''_1}(g) - \omega_1^{\mu'_1}(g) s} e_1.$$
Comme $\Fil^1 \calM_1$ est stable par l'action de $G$, cet élément doit 
être dans $\Fil^1 \calM_1$, ce qui implique que le facteur entre crochets 
est nul (puisque l'on a supposé $s$ de degré strictement inférieur à 
$e$). Cette condition se réécrit $\hat g(s) = \omega_1^{\mu'_1 - 
\mu''_1}(g) s$. Ainsi si $n_2$ est l'unique entier de $[0,p^r-1[$ congru 
à $\mu'_1 - \mu''_1$ modulo $(p^r-1)$ (on rappelle que $e = p^r - 1$), 
on a nécessairement $s = \lambda_2 u^{n_2}$ avec $\lambda_2 \in E$. Donc :

\begin{lemme}
\label{lem:filt}
Si $2 \not\in J$, le sous-module $\Fil^1 \calM_1$ est engendré par les 
deux éléments $e_1$ et $u^e f_1$.

Si $2 \in J$, il existe $\lambda_2 \in E$ tel que $\Fil^1 \calM_1$ 
soit engendré par les deux éléments $u^e e_1$ et $f_1 + \lambda_2 
u^{n_2} e_1$ où $n_2$ est l'unique entier de $[0, p^r-1[$ congru à
$\mu'_1 - \mu''_1$ modulo $(p^r-1)$.
\end{lemme}

\paragraph{Pour les $\mathbf{\calM_i}$, $\mathbf{i > 1}$}

Posons à présent $f_2 = \phi_1(u^e f_1)$ si $2 \not\in J$ et $f_2 = 
\phi_1 (f_1 + \lambda_2 u^{n_2} e_1)$ sinon. Il est clair que $f_2$ se 
réduit sur $\bar f_2$ dans $\calM''$. De plus, en utilisant la 
commutation de $[g]$ et $\phi_1$, on obtient :
$$[g](f_2) = \omega_2^{\mu'_2}(g) f_2$$
(notez que $\omega_i^{\mu'_i}$ ne dépend pas de $i$). En réappliquant
l'argument de la preuve du lemme \ref{lem:filt}, on voit que $\Fil^1 
\calM_2$ est engendré par $e_2$ et $u^e f_2$ si $3 \not\in J$, ou 
sinon par $u^e e_2$ et $f_2 + \lambda_3 u^{n_3} e_1$ avec $\lambda_3 \in 
E$, $0 \leq n_3 < p^r - 1$ et $n_3 \equiv \mu'_2 - \mu''_2 \pmod 
{p^r-1}$.

Ainsi de suite, on construit $f_3$, $f_4$, \emph{etc.}, jusqu'à arriver
à $f_{r+1}$. Malheureusement, rien ne nous dit que celui-ci est égal à
$a f_1$. Cependant l'image de $f_{r+1}$ dans $\calM''$ est bien $a \bar 
f_1$, d'où $f_{r+1} = a f_1 + s e_1$ pour un certain $s \in 
E[u]/u^{ep}$. On doit par ailleurs avoir :
$$[g](f_{r+1}) = \omega_1^{\mu'_1}(g) f_{r+1}$$
d'où l'on déduit que $s = \beta u^{pn}$ ($\beta \in E$) où $n$ est celui 
du lemme \ref{lem:galois}. Ce même lemme nous précise que l'on peut 
modifier $f_1$ en lui ajoutant une quantité de la forme $\alpha u^{pn}$. 
Voyons comment cela modifie $f_2$.

\medskip

Si $2 \not\in J$, il est facile de voir que $f_2$ n'est pas modifié. 
Ainsi en remplaçant $f_1$ par $\frac{f_{r+1}}a$, on peut supposer que 
$\beta = 0$, c'est-à-dire que $f_{r+1} = a f_1$. De même, si $J \neq S$,
on peut reprendre tout le raisonnement précédent en commençant non pas
à $1$ mais à un entier $i_0$ tel que $i_0 + 1 \not\in J$.

Supposons donc $J = S$. Dans ce cas, $\mu'_1 = \sum_{j=1}^r p^{r-j} 
(1-b_{j+1})$ et $\mu''_1 = 0$. On en déduit la valeur de $n_2$ :
$$n_2 = p^r - 1 - \sum_{j=1}^r p^{r-j} (b_{j+1} - 1)$$
ce qui donne en particulier $n_2 \equiv -b_1 \pmod p$. D'autre part, 
$pn \equiv n_2 \pmod{p^r-1}$ (puisqu'ils sont tous les deux congrus à
$\mu'_1 - \mu''_1$). Écrivons $pn = n_2 + q(p^r - 1)$ avec $q \geq 0$. 
Réduisant cette égalité modulo $p$, il reste $q \equiv n_2 \equiv -b_1 
\pmod p$ et donc en particulier $q$ ne peut valoir ni $0$, ni $1$.
Ainsi $q \geq 2$, ce qui implique que $f_2$ n'est pas modifié lorsque
l'on ajoute à $f_1$ une quantité de la forme $\alpha u^{pn}$. On conclut
comme précédemment.

\paragraph{Conclusion}

Au final, la structure de $\calM$ est donnée par la proposition 
suivante :

\begin{prop}
\label{prop:desccalM}
Avec les notations précédentes, il existe, pour tout $i$ une 
$(E[u]/u^{ep})$-base $(e_i, f_i)$ de $\calM_i$ et des éléments 
$\lambda_i \in E$ tels que $\lambda_i = 0$ si $i \not\in J$ et :
\begin{itemize}
\item[i)] $[g](e_i) = \omega_i^{\mu''_i}(g) e_i$ et $[g](f_i) = 
\omega_i^{\mu'_i}(g) f_i$ pour tout $i \in S$ ;
\item[ii)] $\Fil^1 \calM_i$ est engendré par $e'_i = u^{e 
\1_J(i+1)} e_i$ et $f'_i = u^{e \bar \1_J(i+1)} f_i + \lambda_{i+1} 
u^{n_{i+1}} e_i$ où $n_{i+1}$ est le reste de la division euclidienne de 
$\mu'_i - \mu''_i$ par $p^r-1$ ;
\item[iii)] $\phi_1(e'_i) = \underline b (i) e_{i+1}$, $\phi_1(f'_i) = 
\underline a(i) f_{i+1}$ ;
\item[iv)] $N(e_i) = 0$, $N(f_i) = - \frac{\underline b(i-1)} 
{\underline a(i-1)} n_i \lambda_i u^{p n_i} e_i$.
\end{itemize}

De plus, le morphisme $\calM'' \to \calM$ (resp. $\calM \to \calM'$) est 
donné sur cette description par $e_i \mapsto e_i$ (resp. $e_i 
\mapsto 0$, $f_i \mapsto \bar f_i$).
\end{prop}

\begin{proof}
On a déjà tout prouvé sauf la forme de l'opérateur de monodromie. Comme
$\kappa = 2$, il suffit de vérifier que le $N$ de la proposition prolongé 
par la relation de
Leibniz vérifie $N(\calM) \subset u\calM$ et commute à $\phi_1$.
La première condition est équivalente à $n_i > 0$ pour tout $i \in S$,
et en écrivant explicitement la relation de commutation, on montre que
la seconde condition est impliquée par les inégalités $p n_i \geq e$
pour tout $i \in S$.

Nous reportons la preuve de ce deuxième énoncé combinatoire (qui implique
à l'évidence le premier) à une sous-section ultérieure (lemme 
\ref{lem:inegn}).
\end{proof}

\subsection{Interlude : défaut de pleine fidélité}

Oublions un instant le fil de la démonstration pour nous concentrer sur
une question différente. Soient $\calA$ et $\calB$ deux objets de 
$\EBrModdd^{\kappa-1}$ de rang $1$. On note $a_i$, $\alpha_i$ et 
$\alpha_{\f,i}$ (resp. $b_i$, $\beta_i$ et $\beta_{\f,i}$) les
invariants numériques associés à $\calA$ et $\calB$ par la proposition
\ref{prop:rang1}, et $A_i$ (resp. $B_i)$ une base correspondant à ces
invariants. En particulier, tout au long de cet interlude, les $a_i$ et 
les $b_i$ n'ont aucun rapport avec ceux que l'on manipulait jusqu'à 
présent.

\begin{prop}
\label{prop:ordre}
On conserve les notations précédentes et on suppose donné un 
isomorphisme $f : \Tst(\calB) \to \Tst(\calA)$. Alors, il existe un
morphisme non nul (dans la catégorie $\EBrModdd^{\kappa-1}$) $\calA \to 
\calB$ si, et seulement si $\beta_{\f,i} \geq \alpha_{\f,i}$ pour 
tout $i \in S$.

Le cas échéant, il existe $\lambda \in E^\star$ tel que le morphisme 
$\hat f : \calA \to \calB$ défini par $A_i \mapsto \lambda 
u^{\beta_{\f,i} - \alpha_{\f,i}} B_i$ satisfasse $\Tst(\hat f) = 
f$.
\end{prop}

\begin{proof}
Le fait que $\Tst(\calA)$ soit isomorphe à $\Tst(\calB)$ entraîne,
grâce à la proposition \ref{prop:gal1}, d'une part que les éléments de 
$E^\star$ associés à $\calA$ et $\calB$ par la proposition 
\ref{prop:rang1} coïncident, et d'autre part que la congruence
$\alpha_i + \alpha_{\f,i} \equiv \beta_i + \beta_{\f,i} \pmod
{p^r-1}$ est satisfaite.

Soit $\hat f : \calA \to \calB$ un morphisme non nul dans la catégorie
$\EBrModdd^{\kappa-1}$. L'image de $A_i$ est nécessairement de la
forme $s_i B_i$ avec $s_i \in E[u]/u^{ep}$. Si l'un des $s_i$ est nul,
la condition de commutation à $\phi_{\kappa-1}$ implique que $s_{i+1}$ 
l'est également, et par récurrence que $\hat f$ est lui aussi nul. Comme 
ce cas a été exclu, tous les $s_i$ sont non nuls. Notons $v_i$ la \og 
valuation $u$-adique \fg\ de $s_i$ : c'est un entier compris entre $0$
et $ep-1$. En utilisant à nouveau, la commutation à $\phi_{\kappa-1}$,
on obtient les relations $v_{i+1} = p v_i + p a_i - p b_i$. Celles-ci
forment un système qui se résout aisément et conduit à $v_i = 
\beta_{\f,i} - \alpha_{\f,i}$. Ceci démontre le sens direct de
l'équivalence de la proposition.

Pour la réciproque, on considère le morphisme défini $\hat f'$
par $A_i \mapsto u^{\beta_{\f,i} - \alpha_{\f,i}} B_i$. Pour
tout $i$, on a par hypothèse $0 \leq a_i \leq e(\kappa-1) \leq e(p-2)$.
On en déduit :
$$\alpha_{\f,i} \leq \frac{e(p-2)(p^r + p^{r-1} + \cdots + p)}{p^r-1}
= \frac{ep(p-2)}{p-1} < ep$$
pour tout $i$. De même $\beta_{\f,i} < ep$, d'où il suit que la différence
$\beta_{\f,i} - \alpha_{\f,i}$ est elle aussi strictement inférieure
à $ep$. Ainsi le morphisme $\hat f'$ est non nul. Il reste à vérifier
qu'il commute aux structures supplémentaires, mais cela ne pose plus aucune
difficulté particulière. (Pour la donnée de descente, il faut utiliser 
la congruence $\alpha_i + \alpha_{\f,i} \equiv \beta_i + \beta_{\f,i} 
\pmod {p^r-1}$.)

Reste à prouver la dernière assertion. La preuve de la proposition 
\ref{prop:gal1} (dont nous avons seulement expliqué les idées) montre
en fait que $\Tst(\hat f')$ est un isomorphisme. Comme $\Tst(\calA)$ et
$\Tst(\calB)$ sont des $E$-espaces vectoriels de dimension $1$, on a
nécessairement $f = \lambda \Tst(\hat f')$ pour un $\lambda \in 
E^\star$. Il suffit donc de poser $\hat f = \lambda \hat f'$.
\end{proof}

\noindent
{\it Remarque.}
 Quitte à modifier la base $A_i$, on peut supposer dans 
l'énoncé de la proposition que $\lambda = 1$. C'est ce que nous ferons
par la suite.

\bigskip

L'énoncé suivant montre même que dans le \og mauvais \fg\ cas de la 
proposition, la situation n'est pas désespérée :

\begin{prop}
\label{prop:max}
On conserve les notations précédentes et on suppose toujours donné un
isomorphisme $f : \Tst(\calB) \to \Tst(\calA)$. Alors il existe un 
objet $\calC$ de $\EBrModdd^{\kappa-1}$ et des morphismes $\hat f_
\calA : \calA \to \calC$ et $\hat f_\calB : \calB \to \calC$ induisant
des isomorphismes \emph{via} $\Tst$ et tels que $\Tst(\hat f_\calA) 
\circ \Tst(\hat f_\calB)^{-1} = f$.
\end{prop}

\begin{proof}
Évidemment si $\beta_{\f,i} \geq \alpha_{\f,i}$ pour tout $i \in S$, il 
suffit de prendre $\calC = \calB$, $\hat f_\calB = \id$ et d'appliquer
la proposition \ref{prop:ordre} pour obtenir $\hat f_\calA$. De même,
si $\alpha_{\f,i} \geq \beta_{\f,i}$ pour tout $i \in S$, on peut 
conclure avec $\calC = \calA$.

Lorsque les $r$-uplets $(\alpha_{\f,i})$ et $(\beta_{\f,i})$ 
ne sont pas comparables, c'est plus compliqué. L'idée consiste à 
construire un $\calC$ dont les invariants associées $\gamma_{\f,i}$
vérifient $\gamma_{\f,i} = \max (\alpha_{\f,i}, \beta_{\f,i})$. Pour 
y parvenir, on pose $n_i = \frac 1 p \max(\beta_{\f,i} - 
\alpha_{\f,i}, 0)$ et $c_i = a_i + p n_i - n_{i+1}$. Ces nombres sont
toujours des entiers, car on voit directement sur la définition
que $\alpha_{\f,i}$ et $\beta_{\f,i}$ sont des multiples de $p$. 
Montrons que $c_i$ est positif ou nul. Cela ne pose pas de problème si 
$n_{i+1} = 0$. Si par contre $n_{i+1} > 0$, on a $n_{i+1} = \frac 
{\beta_{\f,i+1} - \alpha_{\f,i+1}} p$ puis :
$$a_i + pn_i \geq a_i + \beta_{\f,i} - \alpha_{\f,i} = b_i + \frac 
{\beta_{\f,i+1} - \alpha_{\f,i+1}} p = b_i + n_{i+1} \geq n_{i+1}$$
la première égalité se vérifiant simplement en remplaçant
$\alpha_{\f,i}$, \emph{etc}. par leurs définitions. De même, on 
prouve $c_i \leq e(\kappa-1)$ : si $n_i = 0$, le résulat est 
immédiat et sinon on écrit :
$$a_i + pn_i = a_i + \beta_{\f,i} - \alpha_{\f,i} = b_i + \frac 
{\beta_{\f,i+1} - \alpha_{\f,i+1}} p \leq b_i + n_{i+1} \leq 
e(\kappa-1) + n_{i+1}.$$
(En réalité, on a même obtenu $\min(a_i, b_i) \leq c_i \leq 
\max(a_i,b_i)$.) Un calcul simple
montre par ailleurs que les $\gamma_{\f,i}$ associés aux $c_i$ vérifient 
$\gamma_{\f,i} = \alpha_{\f,i} + pn_i = \max (\alpha_{\f,i}, 
\beta_{\f,i})$, c'est-à-dire la condition rechechée. Définissons
à présent $\gamma_i = \alpha_i + \alpha_{\f,i} - \gamma_{\f,i} \in 
\Z/(p^r-1)\Z$ et $\calC$ l'objet de $\EBrModdd^{\kappa-1}$ associé aux
nombres $c_i$, $\gamma_i$ et à l'élément de $E^\star$ correspondant
à $\calA$ (ou ce qui revient au même, à $\calB$). La proposition 
\ref{prop:gal1} montre que $\Tst(\calA)$ et $\Tst(\calC)$ sont 
isomorphes et il ne reste plus, pour conclure la preuve, qu'à appliquer 
deux fois la proposition \ref{prop:ordre}.
\end{proof}

\noindent
{\it Remarque.}
Les résultats précédents ne doivent pas surprendre le lecteur. Ils sont 
simplement l'illustration en termes de la catégorie $\EBrModdd^{\kappa-1}$ de 
propriétés analogues sur les schémas en groupes de type $(p, \ldots, p)$ 
connues depuis Raynaud (voir \cite{raynaud}, \S 2.2), la relation 
d'ordre (partiel) sur les prolongements d'un $p$-groupe sur $K$ 
correspondant ici à l'ordre produit sur les $\alpha_{\f,i}$.

\subsection{Fin de la démonstration}
\label{subsec:demfin}

Revenons à nos moutons et reprenons en particulier les extensions 
$0 \to M'' \to M \to M' \to 0$ et $0 \to \calM'' \to \calM \to \calM' 
\to 0$ données par les propositions \ref{prop:descM} et 
\ref{prop:desccalM}. Voyons les objets $\calM''$, $\calM$ et $\calM'$
comme appartenant à la catégorie $\EBrModdd^{p-2}$ et notons $\calN''$, 
$\calN$ to $\calN'$ les objets de $\EBrModdd^{p-2}$ associés 
respectivement à $M''$, $M$ et $M'$ par la recette de \ref{subsec:modp}.
Si on définit les constantes :
$$m'_i = e(p-3) + e \bar \1_J(i+1) \quad ; \quad
m''_i = e(p-3) + e \1_J(i+1)$$
$$n'_i = e(p-2-b_i \1_J(i)) \quad ; \quad
n''_i = e(p-2-b_i \bar \1_J(i))$$
les objets $\calM$ et $\calN$ se décrivent explicitement comme suit :
$$\begin{array}{l}
\left\{ \begin{array}{l}
\calM_i = E[u]/u^{ep} \cdot e_i \oplus E[u]/u^{ep} \cdot f_i \\
\Fil^{p-2} \calM_i \text{ est engendré par } u^{m''_i} e_i 
\text{ et } u^{m'_i} f_i + \lambda_{i+1} u^{e(p-3) + n_{i+1}} e_i 
\\
\phi_{p-2}(u^{m''_i} e_i) = \underline b(i) e_{i+1}, \,
\phi_{p-2}(u^{m'_i} f_i + \lambda_{i+1} u^{e(p-3) + n_{i+1}} e_i) 
= \underline a(i) f_{i+1} \\
N(\calE_i) = 0 , \,
N(\calF_i) = -\frac{\underline b (i-1)} {\underline a (i-1)} n_i 
\lambda_i u^{pn_i} e_i \\{}
[g](\calE_i) = \omega_i^{\mu''_i}(g) \calE_i, \,
[g](\calF_i) = \omega_i^{\mu'_i}(g) \calF_i
\end{array} \right. \\ \\
\left\{ \begin{array}{l}
\calN_i = E[u]/u^{ep} \cdot E_i \oplus E[u]/u^{ep} \cdot F_i \\
\Fil^{p-2} \calN_i \text{ est engendré par } u^{n''_i} E_i 
\text{ et } u^{n'_i} F_i \\
\phi_{p-2}(u^{n''_i} E_i) = \underline b(i) E_{i+1}, \,
\phi_{p-2}(u^{n'_i} F_i) = \underline a(i) (F_{i+1} -
\bar \Lambda_{i+1} E_{i+1}) \\
N(E_i) = 0 , \, N(F_i) = 0 \\{}
[g](E_i) = E_i, \, [g](F_i) = F_i
\end{array} \right. \end{array}$$
où $\bar \Lambda_i$ désigne la réduction modulo $p$ de $\Lambda_i$. On
rappelle que lorsque $i \not\in J$, on a $\Lambda_i = 0$ et $\lambda_i
= 0$. On rappelle également que l'objet $\calM''$ (resp. $\calN''$) est
le sous-objet de $\calM$ (resp. $\calN$) engendré par les $e_i$ (resp. 
les $E_i$) et que $\calM'$ (resp. $\calN'$) est le quotient par ce
sous-objet. On note toujours $\bar f_i$ (resp. $\bar F_i$) les images de 
$f_i$ (resp. $F_i$) dans ce quotient.

\paragraph{Stratégie de la preuve}

Les objets précédents viennent avec des identifications $\Tst(\calM') \simeq 
E(\psi') \simeq \Tst(\calN')$ et $\Tst(\calM'') \simeq E(\psi'') 
\simeq \Tst(\calN'')$. Appelons $f'$ et $f''$ les isomorphismes composés 
$\Tst(\calM') \to \Tst(\calN')$ et $\Tst(\calM'') \to \Tst(\calN'')$. 
Pour conclure la preuve qui nous intéresse, il suffit de montrer que si 
$\underline a(i-1) \bar \Lambda_i = \underline b(i-1) \lambda_i$, 
il existe un (iso)morphisme $f : \Tst(\calM) \to \Tst(\calN)$ rendant 
commutatif le diagramme suivant :
$$\xymatrix @C=40pt @R=13pt {
0 \ar[r] & \Tst(\calM') \ar[r] \ar[d]_-{f'}^-{\sim} & \Tst(\calM) \ar[r] 
\ar@{.>}[d]^-{f} & \Tst(\calM'') \ar[r] \ar[d]^-{f''}_-{\sim} & 0 \\
0 \ar[r] & \Tst(\calN') \ar[r] & \Tst(\calN) \ar[r] & \Tst(\calN'') 
\ar[r] & 0 }$$
à lignes exactes. Pour cela, on n'a pas envie de calculer 
explicitement les représentations galoisiennes $\Tst(\calM)$ et 
$\Tst(\calN)$. L'idée, au contraire, est d'essayer de relever le
diagramme précédent au niveau des objets de $\EBrModdd$. On commence 
modestement par relever les isomorphismes $f'$ et $f''$ et c'est là 
qu'interviennent les résultats de l'interlude.

Précisément, d'après la proposition \ref{prop:max}, le morphisme $f'$ 
(resp. $f''$) se relève \emph{via} des morphismes $\hat f'_\calN : 
\calN' \to \calC'$ et $\hat f'_\calM : \calM' \to \calC'$ (resp. 
$\hat f''_\calN : \calN'' \to \calC''$ et $\hat f''_\calM : \calM'' \to 
\calC''$) pour un certain objet $\calC'$ (resp. $\calC''$) de la 
catégorie $\EBrModdd^{p-2}$. On cherche donc au final à construire un 
objet $\calC \in \EBrModdd^{p-2}$ qui s'insère dans le diagramme 
commutatif suivant :
\begin{equation}
\label{eq:diagramme}
\raisebox{0.5\depth}{\xymatrix @C=40pt @R=16pt {
0 \ar[r] & \calN'' \ar[r] \ar[d]_-{\hat f''_\calN} & \calN \ar[r] 
\ar@{.>}[d]^-{\hat f_\calN} & \calN' \ar[r] \ar[d]^-{\hat f'_\calN} 
& 0 \\
0 \ar[r] & \calC'' \ar@{.>}[r] & \calC \ar@{.>}[r] & \calC' \ar[r] 
& 0 \\
0 \ar[r] & \calM'' \ar[r] \ar[u]^-{\hat f''_\calM} & \calM \ar[r] 
\ar@{.>}[u]_-{\hat f_\calM} & \calM' \ar[r] \ar[u]_-{\hat f'_\calM} 
& 0 }}
\end{equation}
où les flèches en pointillés restent aussi à définir. Remarquons tout de
suite que si l'on y parvient, on aura terminé la démonstration puisque
le lemme des cinq prouvera immédiatement que les applications $\Tst(\hat
f_\calM)$ et $\Tst(\hat f_\calN)$ sont des isomorphismes.

\paragraph{Description de $\mathbf \calC'$ et $\mathbf \calC''$}

Notons $\mu'_{\f,i}$, $\mu''_{\f,i}$, $\nu'_{\f,i}$ et $\nu''_{\f,i}$ 
les sommes correspondantes respectivement aux nombres $m'_i$, $m''_i$, 
$n'_i$ et $n''_i$ par la formule de la proposition \ref{prop:rang1}.
Elle valent :
$$\mu'_{\f,i} = p(p-3) + \sum_{j=0}^{r-1} p^{r-j} \bar \1_J(i+j+1)
\quad ; \quad
\mu''_{\f,i} = p(p-3) + \sum_{j=0}^{r-1} p^{r-j} \1_J(i+j+1)$$
$$\nu'_{\f,i} = \sum_{j=0}^{r-1} p^{r-j} (p-2-b_{i+j} \1_J(i+j))
\quad ; \quad
\nu''_{\f,i} = \sum_{j=0}^{r-1} p^{r-j} (p-2 - b_{i+j} \bar 
\1_J(i+j)).$$

\begin{lemme}
\label{lem:ineg}
On a $i \in J$ si, et seulement si $\mu'_{\f,i} > \nu'_{\f,i}$ si, et 
seulement si $\mu''_{\f,i} \leq \nu''_{\f,i}$.
\end{lemme}

\begin{proof}
Des formules donnant les valeurs de $\mu'_{\f,i}$ et $\nu'_{\f,i}$, on
déduit l'expression suivante :
$$\mu'_{\f,i} - \nu'_{\f,i} = \sum_{j=0}^{r-1} p^{r-j} \big[ 
b_{i+j} \1_J(i+j) - \1_J(i+j+1) \big].$$
Comme $2 \leq b_{i+j} \leq p-2$, le facteur entre crochets a toujours 
une valeur absolue inférieure ou égale à $p-1$. Ainsi, le 
signe de la somme est donné par le premier terme (\emph{i.e.} celui pour 
lequel $j$ est le plus petit) non nul de celle-ci. Si $i \in J$, ce 
terme est celui obtenu pour $j=0$, et la somme est donc 
strictement positive, \emph{i.e.} $\mu'_{\f,i} > \nu'_{\f,i}$. Supposons 
maintenant $i \not\in J$. Si $J = \emptyset$, alors la somme est nulle 
et $\mu'_{\f,i} = \nu'_{\f,i}$. Sinon, le premier terme non nul 
correspond au plus petit $j$ tel que $i+j+1 \in J$ et le facteur entre 
crochets associé est négatif (puisque $i+j \not\in J$ et donc le terme 
avec $b_{i+j}$ n'apparaît pas). L'équivalence s'ensuit.

L'inégalité avec les $\mu''_{\f,i}$ et $\nu''_{\f,i}$ se traite de même.
\end{proof}

Posons pour simplifier les écritures qui vont suivre $s_i = 
\frac{\mu'_{\f,i} - \nu'_{\f,i}} p$ et $r_i = \frac{\mu''_{\f,i} - 
\nu''_{\f,i}} p$. À partir du lemme précédent, une étude de la preuve
de la proposition \ref{prop:max} montre que $\calC'$ et $\calC''$ se
décrivent comme suit :
$$\left\{\begin{array}{l}
\calC'_i = E[u]/u^{ep} \cdot \bar \calF_i \\
\Fil^1 \calC'_i = u^{c'_i} \calC'_i \\
\phi_1(u^{c'_i} \bar \calF_i) = \underline a (i) \bar \calF_{i+1} 
\\{}
[g] (\bar \calF_i) = \omega_i^{\gamma'_i}(g) \bar \calF_i
\end{array}\right.
\qquad 
\left\{\begin{array}{l}
\calC''_i = E[u]/u^{ep} \cdot \calE_i \\
\Fil^1 \calC''_i = u^{c''_i} \calC''_i \\
\phi_1(u^{c''_i} \calE_i) = \underline b (i) \calE_{i+1} \\{}
[g] (\calE_i) = \omega_i^{\gamma''_i}(g) \calE_i
\end{array}\right.$$
avec :
$$\gamma'_i = p s_i \1_J(i) \quad ; \quad \gamma''_i = - p r_i 
\bar \1_J(i)$$
$$c'_i = - s_{i+1} \1_J(i+1) + ps_i \1_J(i) + e(p-2-b_i \1_J(i))$$
$$c''_i = - r_{i+1} \bar \1_J(i+1) + pr_i \bar \1_J(i) + e(p-2-b_i 
\bar \1_J(i))$$

\paragraph{Quelques lemmes préparatoires}

Avant de donner la construction de $\calC$, nous regroupons dans ce
paragraphe les identités combinatoires qu'il est bon de garder en
tête pour faire les vérifications à venir :

\begin{lemme}
\label{lem:ident}
On a les identités suivantes :
\begin{itemize}
\item[$\bullet$] $p r_i - r_{i+1} = e \big[ b_i \bar \1_J(i) - \bar 
\1_J(i+1) \big]$
\item[$\bullet$] $p s_i - s_{i+1} = e \big[ b_i \1_J(i) - \1_J(i+1) 
\big]$ 
\item[$\bullet$] $c'_i = s_{i+1} \bar \1_J(i+1) - ps_i \bar \1_J(i) + 
e(p-3) + e \bar \1_J(i+1)$
\item[$\bullet$] $c''_i = r_{i+1} \1_J(i+1) - pr_i \1_J(i) + e(p-3) + e 
\1_J(i+1)$
\item[$\bullet$] $n_i = r_i - s_i + e \1_J(i)$ 
\end{itemize}
où on rappelle que $n_i$ est défini dans la proposition 
\ref{prop:desccalM}. Si de plus $i+1 \in J$, on a aussi :
\begin{itemize}
\item[$\bullet$] $c''_i - s_{i+1} = e(p-3) - pr_i \1_J(i) + n_{i+1}$
\item[$\bullet$] $c''_i - s_{i+1} - c'_i = p s_i \bar \1_J(i) - p r_i
\1_J(i) + n_{i+1} \equiv n_{i+1} \pmod p$
\end{itemize}
\end{lemme}

\begin{proof}
Les deux premières égalités résultent directement des définitions, les
deux suivantes en sont des conséquences immédiates, et les deux dernières
s'obtiennent facilement à partir des quatre premières. Reste donc simplement 
à montrer la cinquième. Pour cela, on remarque que l'on a les congruences $r_i 
\equiv p r_{i-1} \equiv - \mu''_{i-1} \pmod{p^r-1}$ et $s_i \equiv p 
s_{i-1} \equiv - \mu'_{i-1} \pmod{p^r-1}$. Ainsi $n_i \equiv r_i - s_i 
\pmod{p^r-1}$. Grâce à l'écriture :
$$r_i - s_i = \sum_{j=0}^{r-1} p^{r-j-1} \big[ b_{i+j} (\bar \1_J - 
\1_J) (i+j) - (\bar \1_J - \1_J) (i+j+1)\big]$$
on s'aperçoit que $|r_i - s_i| < e$.
Supposons $i \in J$. D'après le lemme \ref{lem:ineg}, on a $r_i \leq 0$ 
et $s_i > 0$. Ainsi $r_i - s_i < 0$ et $n_i = r_i - s_i + e$. Si au 
contraire, $i \not\in J$, on a $r_i - s_i \geq 0$ et puis $n_i = r_i - 
s_i$.
\end{proof}

\begin{lemme}
\label{lem:ineg2}
Si $i+1 \in J$, alors $c''_i \geq s_{i+1}$ et $c'_i + s_{i+1} \leq 
e(p-2)$.
\end{lemme}

\begin{proof}
Par le lemme \ref{lem:ident}, la différence $c''_i - s_{i+1}$ 
s'exprime comme une somme de trois termes positifs. La première 
inégalité s'ensuit. Pour la seconde, l'expression :
$$c'_i + s_{i+1} = e(p-2) + (p s_i - e b_i) \1_J(i).$$
montre qu'il suffit donc de prouver $p s_i - e b_i \leq 0$
pour $i \in J$. Or, dans ce cas :
$$p s_i - e b_i = - p^r + b_i + \sum_{j=1}^{r-1} p^{r-j} \big[ b_{i+j} 
\1_J(i+j) - \1_J(i+j+1) \big].$$
La conclusion s'obtient en remarquant que $b_i$ ainsi que tous les 
facteurs entre crochets sont majorés en valeur absolue par $p-1$.
\end{proof}

\begin{lemme}
\label{lem:inegn}
Pour tout $i \in S$, on a $p n_i \geq e$.
\end{lemme}

\begin{proof}
Les égalités du lemme \ref{lem:ident} impliquent la congruence :
$$n_{i+1} = r_{i+1} - s_{i+1} + e \1_J(i+1) \equiv b_i(\bar \1_J - \1_J)(i) 
- \bar \1_J (i+1) \pmod p.$$
De $2 \leq b_i \leq p-2$, on déduit que $n_{i+1}$ n'est jamais un multiple 
de $p$. Par ailleurs, par le même lemme, on obtient les congruences
$r_{i+1} \equiv p r_i \pmod e$, $s_{i+1} \equiv p s_i \pmod e$ et
$n_i \equiv r_i - s_i \pmod e$. Il s'ensuit :
\begin{equation}
\label{eq:congr}
n_{i+1} \equiv p n_i \pmod e
\end{equation}
Fixons $i \in S$ et supposons par l'absurde que $p n_i < e$. Comme, par
hypothèse, on a $n_i \geq 0$ et $0 \leq n_{i+1} < e$, la congruence
(\ref{eq:congr}) entraîne l'égalite $p n_i = n_{i+1}$, ce qui contredit
le fait que $p$ ne divise pas $n_{i+1}$.
\end{proof}

\paragraph{Définition de $\mathbf \calC$}

On est enfin prêt à donner la définition du fameux objet $\calC$. Il est 
obtenu par les formules suivantes :
$$\left\{ \begin{array}{l}
\calC_i = E[u]/u^{ep} \cdot \calE_i \oplus E[u]/u^{ep} \cdot \calF_i \\
\Fil^{p-2} \calC_i \text{ est engendré par } u^{c''_i} \calE_i 
\text{ et } u^{c'_i} \calF_i + \lambda_{i+1} u^{c''_i - s_{i+1}} \calE_i 
\\
\phi_{p-2}(u^{c''_i} \calE_i) = \underline b(i) \calE_{i+1}, \,
\phi_{p-2}(u^{c'_i} \calF_i + \lambda_{i+1} u^{c''_i - s_{i+1}} \calE_i) 
= \underline a(i) \calF_{i+1} \\
N(\calE_i) = 0 , \,
N(\calF_i) = -\frac{\underline b (i-1)} {\underline a (i-1)} n_i 
\lambda_i u^{pn_i - pr_i} \calE_i \\{}
[g](\calE_i) = \omega_i^{\gamma''_i}(g) \calE_i, \,
[g](\calF_i) = \omega_i^{\gamma'_i}(g) \calF_i
\end{array} \right.$$
La première inégalité du lemme \ref{lem:ineg2} assure que le terme 
$\lambda_{i+1} u^{c''_i - s_{i+1}}$ qui apparaît dans l'écriture précédente 
a bien un sens (on rappelle que $\lambda_i$ est nul pour $i \not\in J$), 
alors que la deuxième inégalité du même lemme implique l'inclusion
$u^{e(p-2)} \calC \subset \Fil^{p-2} \calC$. Les autres axiomes de la
définition de la catégorie $\EBrModdd^{p-2}$ se vérifient sans grande
difficulté (éventuellement en utilisant les égalités et la congruence
du lemme \ref{lem:ident}). Bien évidemment 
le morphisme $\calC''_i \to \calC_i$ est donné par $\calE_i \mapsto 
\calE_i$ alors que celui $\calC_i \to \calC'_i$ est donné par $\calE_i 
\mapsto 0$, $\calF_i \mapsto \bar \calF_i$.

\bigskip

Voici enfin le lemme qui conclut la preuve :

\begin{lemme}
Les applications $E[u]/u^{ep}$-linéaires $f_\calM : \calM \to \calC$ et 
$f_\calN : \calN \to \calC$ définis par $f_\calM(e_i) = u^{-pr_i 
\1_J(i)} \calE_i$, $f_\calM(f_i) = u^{-ps_i \bar \1_J(i)} \calF_i$, 
$f_\calN(E_i) = u^{pr_i \bar \1_J(i)} \calE_i$ et $f_\calN(F_i) = 
u^{ps_i \1_J(i)} \calF_i$ définissent des morphismes dans
$\EBrModdd^{p-2}$ qui font commuter le diagramme (\ref{eq:diagramme}).
\end{lemme}

\begin{proof}
C'est une vérification un peu longue mais très facile avec les identités 
regroupées dans le lemme \ref{lem:ident}. (Il est important aussi de 
garder à l'esprit que $\lambda_i$ est nul pour $i \not\in J$ et que
l'on a supposé $\underline a(i-1) \bar \Lambda_i = \underline b(i-1) 
\lambda_i$.) Par voie de conséquence, nous la laissons au lecteur.
\end{proof}


\end{document}